\documentclass[12pt]{article}
\usepackage{makeidx}
\usepackage{amssymb}
\usepackage{amsmath,amsthm}
\usepackage{amsfonts}

\numberwithin{equation}{section}

\begin{document}

\title{Limiting Laws of Linear Eigenvalue Statistics \\
for Unitary Invariant Matrix Models}
\author{L. Pastur \thanks{%
Electronic mail: lpastur@ilt.kharkov.ua} \\
Department of Mathematics, University of Wales Swansea \\
Singleton Park Swansea SA2 8PP UK \thanks{%
Permanent address: Institute for Low Temperatures, 47 Lenin's Ave, 61103
Kharkiv, Ukraine}}
\date{}
\maketitle

\begin{abstract}
\noindent We study the variance and the Laplace transform of the probability
law of linear eigenvalue statistics of unitary invariant Matrix Models of $%
n\times n$ Hermitian matrices as $n\rightarrow \infty $. Assuming that the
test function of statistics is smooth enough and using the asymptotic
formulas by Deift et al for orthogonal polynomials with varying weights, we
show first that if the support of the Density of States of the model
consists of $q\geq 2$ intervals, then in the global regime the variance of
statistics is a quasiperiodic function of $n$ as $n\rightarrow \infty $
generically in the potential, determining the model. We show next that the
exponent of the Laplace transform of the probability law is not in general $%
1/2\ \times $ variance, as it should be if the Central Limit Theorem would
be valid, and we find the asymptotic form of the Laplace transform of the
probability law in certain cases.
\end{abstract}


\newpage

\section{Introduction}

\noindent Random matrix theory deals mostly with eigenvalue distributions of
various ensembles of $n\times n$ random matrices as $n\rightarrow \infty $.
Many questions of this branch of the theory can be formulated in terms of
the Eigenvalue Counting Measure $\mathcal{N}_{n}$, $\mathcal{N}_{n}(\Delta
_{n})$ is the number of eigenvalues in a given part $\Delta _{n}$ of the
spectrum of a matrix $M_{n}$. In a particular case of $n\times n$ Hermitian
matrices eigenvalues are real and we can write for $\Delta _{n}\subset
\mathbb{R}$:
\begin{equation}
\mathcal{N}_{n}(\Delta _{n}):=\sharp \{\lambda _{l}^{(n)}\in \Delta
_{n},\;l=1,...,n\}=\sum_{l=1}^{n}\chi _{\Delta _{n}}(\lambda _{l}^{(n)}),
\label{ECM}
\end{equation}%
where $\chi _{\Delta }$ is the indicator of $\Delta $. A more general object
is a linear eigenvalue statistic, defined as%
\begin{equation}
\mathcal{N}_{n}[\varphi _{n}]:=\sum_{l=1}^{n}\varphi _{n}(\lambda
_{l}^{(n)})=\mathrm{Tr\;}\varphi _{n}(M_{n})=\int_{\mathbb{R}}\varphi
_{n}(\lambda )\mathcal{N}_{n}(d\lambda )  \label{LS}
\end{equation}%
for a certain test function%
\begin{equation}
\varphi _{n}:\mathbb{R}\rightarrow \mathbb{C}  \label{phi_n}
\end{equation}%
It is known that in many cases there exists a scaling of matrix entries
(i.e. a choice of the scale of the spectral axis), such that for a
sufficiently big class of $n$-independent intervals in (\ref{ECM}) (test
functions in (\ref{LS})) the Normalized Counting Measure of eigenvalues
\begin{equation}
N_{n}=\mathcal{N}_{n}\left/ n\right.  \label{NCM}
\end{equation}%
converges weakly to a non-random measure $N$, known as the Integrated
Density of States measure (IDS) of the ensemble%
\begin{equation}
N_{n}\rightarrow N.  \label{IDS}
\end{equation}%
The corresponding scale (asymptotic regime) is called the \textit{global}
(or \textit{macroscopic}). The convergence is either in probability or even
with probability 1. We refer the reader to works \cite{Gi:01,BPS:95,Pa:00},
where this fact is proved and discussed for two most widely studied classes
of random matrix ensembles, the Wigner Ensembles (independent or weakly
dependent entries) and the Matrix Models (invariant matrix probability
laws). Analogous facts are also known for many other ensembles.

Besides, in many cases the measure $N$ possesses a bounded and continuous
density \cite{De-Co:98,Pa:96}:%
\begin{equation}
N(d\lambda )=\rho (\lambda )d\lambda ,  \label{DOS}
\end{equation}%
We call $\rho $ the Density of States (DOS) of ensemble.

These results can be viewed as analogs of the Law of Large Numbers of
probability theory. Hence, a natural next question, important also for
applications, concerns the limiting probability law of fluctuation of the
Normalized Counting Measures of eigenvalues or their linear statistics, i.e.
an analog of the Central Limit Theorem of the probability theory. The
question is not going to be completely trivial, because eigenvalues of
random matrices are strongly dependent, and, as a result, the variance of
the linear statistics (\ref{LS}) with a $C^{1}$ test function does not grow
with $n$ ( see e.g. formulas (\ref{VCr}) -- (\ref{Vbou}) below).
Nevertheless, it was found that in a variety of cases fluctuations of
various spectral characteristics of eigenvalues of random matrix ensembles
are asymptotically Gaussian ( see e.g. \cite%
{Ba-Si:04,Ba:97,Ca:01,Co-Le:95,Di-Ev:01,Gi:01,Gu:05,Jo:98,
Ke-Na:00,KKP,Si-So:97,So:98,So:00,So:02,Sp:87}).In particular, for the
global scale and a $C^{1}$ test function in (\ref{LS}) this requires,
roughly speaking, the same order of magnitude of the entries for the Wigner
Ensembles (see \cite{Gi:01,KKP} for exact conditions, similar to the
Lindeberg condition of the probability theory), and for the Matrix Models
one needs to assume that the support of the IDS (\ref{IDS}) is a connected
interval of the spectral axis \cite{Jo:98}: $\mathrm{supp\;}N=[a,b]$. The
last result was obtained by using variational techniques, introduced in the
random matrix theory in paper \cite{BPS:95} in order to prove (\ref{IDS}).

Being applicable to Matrix Models of all three symmetry classes of the
random matrix theory (real symmetric, Hermitian and quaternion real
matrices), the variational techniques were efficient so far in the study of
fluctuations of eigenvalue statistics only in the case where $\mathrm{supp}%
N=[a,b]$. In this paper we consider only the Matrix Models of Hermitian
matrices, but for a general case of a multi-interval support of $N$:
\begin{equation}
\sigma :=\mathrm{supp\;}N=\bigcup_{l=1}^{q}\;[a_{l},b_{l}],\;q\geq 1.
\label{sig}
\end{equation}%
In this case we can use recent powerful results by Deift et al \cite%
{De-Co:99} on asymptotics of special class of orthogonal polynomials. We
find that if $\varphi $ is real analytic, then the traditional Central Limit
Theorem is not always valid for the case $q\geq 2$. In particular, the
variance and the probability law oscillate in $n$ as $n\rightarrow \infty $,
hence their limiting form depend on a sequence $n_{j}\rightarrow \infty $.
Moreover, the limiting probability laws are not always Gaussian.

The paper is organized as follows. In Section 2 we study the variance of
linear eigenvalue statistics in the global regime, i.e., for $n$-independent
$\varphi $ in (\ref{LS}), confining ourselves mostly to the case of $C^{1}$
test functions $\varphi $. We find that the variance is quasi-periodic in $n$
in general and its frequency module is determined by the charges%
\begin{equation}
\beta _{l}=N([a_{l+1},\infty )),\;l=1,...,q-1,  \label{bl}
\end{equation}%
determined by the IDS of (\ref{IDS}) and its support (\ref{sig}). Hence, the
variance has no limit as $n\rightarrow \infty $ for $q\geq 2$ and its
asymptotic forms are indexed by points of the subset $\mathbb{H}%
^{q-1}\subset \mathbb{T}^{q-1}$, which is the closure of limit points of the
vectors
\begin{equation}
(\{\beta _{1}n\},...,\{\beta _{q-1}n\})\in \mathbb{T}^{q-1},  \label{bn}
\end{equation}%
where $\{\beta _{l}n\},\;l=1,...,q-1$ are the fractional parts of $\beta
_{l}n,\;l=1,...,q-1$. This phenomenon has been already found in certain
cases \cite{APS:01,Da-Co:00}, but we give its general description.

In Section 3 we study the Laplace transform of the probability law of linear
eigenvalue statistics (\ref{LS}) in the global regime, passing to the limit
along a subsequence
\begin{equation}
(\{\beta _{1}n_{j}(x)\},...,\{\beta _{q-1}n_{j}(x)\})\rightarrow x\in
\mathbb{H}^{q-1},  \label{bnl}
\end{equation}%
and confining ourselves to real analytic test functions. We give first a
general formula for the corresponding limit. Since the formula is rather
complex, we consider several particular cases, where we show that the
exponent of the limiting (in the sense (\ref{bnl})) Laplace transform is not
quadratic in $\varphi $, hence the limiting law is not Gaussian (see
formulas (\ref{chf}), (\ref{Flin}), and (\ref{nGau}) -- (\ref{Ax})). This
has to be compared with results of paper \cite{Ch-La:98}, according to which
the limits of variance and the probability law are the same for all
sequences $n_{j}\rightarrow \infty $ (i.e., exist), and the limiting
probability law is Gaussian.

The random matrix theory deals also with two more asymptotic regimes in
addition to the global one. Namely, if, having fixed the global scale
allowing us to prove (\ref{IDS}), we set in (\ref{LS})
\begin{equation}
\varphi _{n}(\lambda )=\varphi \left( (\lambda -\lambda _{0})n^{\alpha
}\right) ,\;0<\alpha <1,\;\lambda _{0}\in \mathrm{supp\;}N,  \label{rint}
\end{equation}%
where $\varphi $ is $n$-independent, we obtain the \textit{intermediate}
regime, and if
\begin{equation}
\varphi _{n}(\lambda )=\varphi \left( (\lambda -\lambda _{0})n\right)
,\;\lambda _{0}\in \mathrm{supp\;}N  \label{rloc}
\end{equation}%
then we have the \textit{local} (or \textit{microscopic}) regime. In Section
4 we discuss the form of variance and the validity of the CLT in these
regimes.

In Appendix A.1 we find the asymptotic form of the variance of the
Eigenvalue Counting Measure (\ref{ECM}) (corresponding to a piece-wise
constant $\varphi _{n}=\chi _{\Delta }$ in (\ref{LS})) for the Gaussian
Unitary Ensemble (GUE). The variance proves to be $O(\log n)$ instead of
being $O(1)$ in the case of a $C^{1}$ test function. This shows that the
Lipschitz condition (\ref{Ho}) on test functions that we use to prove that
the variance is bounded in $n$ is pertinent (for a more precise condition
see \cite{Di-Ev:01}). In Appendix A.2 we compute the variational derivative
of $\beta _{1}$ in the case $q=2$, which is used in Section 4, and discuss
related topics.

We note that a completely rigorous derivation of the results of this paper,
especially those of Section 3, requires rather technical and tedious
arguments. They will not be presented below. Rather, we confine ourselves to
the presentation of results, their discussion, and outline of corresponding
proofs.


\section{Variance of linear eigenvalue statistics}

\subsection{Generalities}

\noindent Recall that unitary invariant Matrix Models are $n\times n$
Hermitian random matrices, defined by the probability law%
\begin{equation}
\mathbf{P}_{n}(dM_{n})=Z_{n}^{-1}\exp \{-n\mathrm{Tr}V(M_{n})\}dM_{n},
\label{PdM}
\end{equation}%
where $M_{n}=\{M_{jk}\}_{j,k=1}^{n},\;M_{jk}=\overline{M_{kj}}$,
\begin{equation*}
dM_{n}=\prod\limits_{j=1}^{n}dM_{jj}\prod\limits_{1\leq j<k\leq n}d\Re
M_{jk}d\Im M_{jk},
\end{equation*}%
and $V:\mathbb{R}\rightarrow \mathbb{R_{+}}$ is a measurable function,
called the potential and such that
\begin{equation}
V:\mathbb{R}\rightarrow \mathbb{R_{+}},\ V(\lambda )\geq (2+\delta )\log
|\lambda |,\ |\lambda |>L,  \label{Vlog}
\end{equation}%
for some positive $\delta $ and $L$.

The limit (\ref{IDS}) can be described as follows \cite{BPS:95,Jo:98}.
Consider the functional:
\begin{equation}
\mathcal{E}_{V}[m]=-\int_{\mathbb{R}}\int_{\mathbb{R}}\log |\lambda -\mu
|m(d\lambda )m(d\mu )+\int_{\mathbb{R}}V(\lambda )m(d\lambda ),  \label{EV}
\end{equation}%
where $m$ is a non-negative unit measure. 

The variational problem, defined by (\ref{EV}), goes back to Gauss and is
called the minimum energy problem in the external field\textit{\ }$V$. The
unit measure $N$ minimizing (\ref{EV}) is called the equilibrium measure in
the external field $V$ because of its evident electrostatic interpretation
as the equilibrium distribution of linear charges on the ideal conductor
occupying the axis $\mathbb{R}$ and confined by the external electric field
of potential $V$. We stress that the respective minimizing procedure
determines both the support $\sigma $ of the measure and the form of the
measure. This should be compared with the variational problem of the theory
of logarithmic potential, where the external field is absent but the support
$\sigma $ is given (see (\ref{Epot})). The minimum energy problem in the
external field (\ref{EV}) arises in various domains of analysis and its
applications (see recent book \cite{Sa-To:97} for a rather complete account
of results and references concerning the problem).

The measure $N$ and its support $\sigma $ are uniquely determined by the
Euler-Lagrange equation of the variational problem \cite{BPS:95,Sa-To:97}:%
\begin{eqnarray}
V_{eff}(\lambda ) &=&F,\;\lambda \in \sigma ,  \label{Vs} \\
V_{eff}(\lambda ) &\geq &F,\;\lambda \notin \sigma ,  \label{Vns}
\end{eqnarray}%
where%
\begin{equation}
V_{eff}(\lambda )=V(\lambda )-2\int_{\sigma }\log |\lambda -\mu |N(d\mu ),
\label{Vef}
\end{equation}%
and $F$ is a constant (the Lagrange multiplier of the normalization
condition $N(\mathbb{R})=1$).

According to \cite{BPS:95} (see also \cite{Jo:98}), if the potential $V$ in (%
\ref{PdM}) -- (\ref{Vlog}) satisfies the local Lipshitz condition%
\begin{equation}
\left\vert V(\lambda _{1})-V(\lambda _{2})\right\vert \leq C|\lambda
_{1}-\lambda _{2}|^{\gamma },\;|\lambda _{1}|,|\lambda _{2}|\leq L,
\label{Lipl}
\end{equation}%
valid for any $L>0$ and some positive $C$ and $\gamma $, then (\ref{IDS})
holds with probability 1, and $N$ is the minimizer of (\ref{EV}). Moreover,
if $V^{\prime }$ satisfies the local Lipshitz condition, and the support (%
\ref{sig}) is a finite union of disjoint finite intervals, then (\ref{DOS})
is valid \cite{De-Co:98,Pa:96} and the Density of States can be written as
\begin{equation}
\rho (\lambda )=P(\lambda )\sqrt{R_{q}(\lambda )},\;\lambda \in \sigma ,
\label{rho}
\end{equation}%
where $P(\lambda )$ is a continuous function,
\begin{equation}
\sqrt{R_{q}(\lambda )}=\left. \sqrt{R_{q}(z)}\right\vert _{z=\lambda
+i0},\;R_{q}(z)=\prod\limits_{l=1}^{q}(z-a_{l})(z-b_{l}),  \label{Xl}
\end{equation}%
and $\sqrt{R_{q}(z)}$ is the branch, determined by the condition: $\sqrt{%
R_{q}(z)}=z^{q}+O(z^{q-1}),\;z\rightarrow \infty $. To obtain these
formulas, provided that the support (\ref{sig}) is given, we differentiate (%
\ref{Vs}), (\ref{Vef}) and obtain the singular integral equation%
\begin{equation}
\mathrm{v.p.}\int_{\sigma }\frac{\rho (\mu )d\mu }{\mu -\lambda }=-\frac{%
V^{\prime }(\lambda )}{2},\;\lambda \in \sigma .  \label{inteq}
\end{equation}%
Then the bounded solution of the equation has the form (\ref{rho}) (see e.g.
\cite{Mu}) in which%
\begin{equation}
P(\lambda )=\frac{1}{2\pi ^{2}}\int_{\sigma }\frac{V^{\prime }(\mu
)-V^{\prime }(\lambda )}{\mu -\lambda }\frac{d\mu }{\sqrt{R_{q}(\mu )}}.
\label{P}
\end{equation}%
The endpoints of the support are rather complex functionals of the potential
in general. Thus, it is of interest to mention a simple case \cite{Bu-Pa:02}.

Let $v:\mathbb{R}\rightarrow \mathbb{R}$ be a monic polynomial of degree $q$
with real coefficients. Assume that for some $g>0$ all zeros of $v^{2}-4g$
are real and simple and set
\begin{equation}
V(\lambda )=\frac{v^{2}(\lambda )}{2gq}.  \label{VBP}
\end{equation}%
Then the DOS of the matrix model (\ref{PdM}) with this potential is
\begin{equation}
\rho (\lambda )=\frac{|v^{\prime }(\lambda )|}{2\pi gq}\left\vert
v^{2}(\lambda )-4g\right\vert ^{1/2},\;\lambda \in \sigma ,  \label{rBP}
\end{equation}%
where
\begin{equation}
\sigma =\{\lambda \in \mathbb{R}:v^{2}(\lambda )\leq 4g\}.  \label{sBP}
\end{equation}%
Besides, in this case we have for the charges (\ref{bl})
\begin{equation}
\beta _{l}=\frac{q-l}{q},\;l=1,...,q-1,  \label{bBP}
\end{equation}%
hence the set $\mathbb{H}^{q-1}$ is $\{0,1/q,...,(q-1)/q\}$.

The case $q=1$ corresponds to the Gaussian Unitary Ensemble and (\ref{rBP})
yields the semi-circle law by Wigner:%
\begin{equation}
V=\frac{\lambda ^{2}}{2g},\;\sigma =[-2g,2g],\;\rho (\lambda )=\frac{1}{2\pi
g}\left\{
\begin{array}{cc}
\sqrt{4g-\lambda ^{2}}, & \lambda \in \sigma , \\
0, & \lambda \notin \sigma .%
\end{array}%
\right.  \label{GUEg}
\end{equation}%
In the case $q=2$ and
\begin{equation}
v(\lambda )=\lambda ^{2}-m^{2},\;m^{2}>2\sqrt{g},  \label{vq2}
\end{equation}%
we have:%
\begin{equation}
\sigma =[-b,-a]\bigcup \;[a,b],\;a=\sqrt{m^{2}-2\sqrt{g}},\;b=\sqrt{m^{2}+2%
\sqrt{g}},  \label{sq2}
\end{equation}%
and
\begin{equation}
\rho (\lambda )=\frac{|\lambda |}{2\pi g}\left\{
\begin{array}{cc}
\sqrt{(b^{2}-\lambda ^{2})(\lambda ^{2}-a^{2})}, & \lambda \in \sigma , \\
0, & \lambda \notin \sigma .%
\end{array}%
\right.  \label{rq2}
\end{equation}

\medskip \noindent

\noindent We will use in this paper the expressions for the variance of
linear statistics (\ref{LS}) and for the Laplace transform of their
probability law via special orthogonal polynomials. The technique dates back
to works by Dyson, Gaudin, Mehta, and Wigner of the 60s (see e.g. \cite%
{Me:92}). Namely, we have for the joint probability density of eigenvalues
of ensemble (\ref{PdM}):
\begin{equation}
p_{n}(\lambda _{1},...,\lambda _{n})=\Big(\det \big\{\psi
_{j-1}^{(n)}(\lambda _{k})\big\}_{j,k=1}^{n}\Big)^{2}\Big/n!,  \label{p_n}
\end{equation}%
where
\begin{equation}
\psi _{l}^{(n)}=e^{-nV/2}P_{l}^{(n)},  \label{psi}
\end{equation}%
and
\begin{equation}
\{P_{l}^{(n)}\}_{l\geq 0}  \label{pl}
\end{equation}%
is the system of orthonormal polynomials with respect to the weight
\begin{equation}
w_{n}=e^{-nV},  \label{wn}
\end{equation}%
so that%
\begin{equation}
\int_{\mathbb{R}}e^{-nV(\lambda )}P_{l}^{(n)}(\lambda )P_{m}^{(n)}(\lambda
)d\lambda =\delta _{l,m},\;l,m=0,1,...  \label{orth}
\end{equation}%
The polynomials satisfy the three - term recurrence relation for $l=0,1,...$%
:
\begin{equation}
r_{l}^{(n)}\psi _{l+1}^{(n)}(\lambda )+s_{l}^{(n)}\psi _{l}^{(n)}(\lambda
)+r_{l-1}^{(n)}\psi _{l-1}^{(n)}(\lambda )=\lambda \psi _{l}^{(n)}(\lambda
),\;r_{-1}^{(n)}=0,  \label{rec}
\end{equation}%
thereby determining a semi - infinite Jacobi matrix:
\begin{equation}
J_{j,k}^{(n)}=r_{j}^{(n)}\delta _{j+1,k}+s_{k}^{(n)}\delta
_{j,k}+r_{j-1}^{(n)}\delta _{j-1,k}.\;j,k=0,1,...  \label{J_n}
\end{equation}%
By using (\ref{p_n}) it can be shown that \cite{Pa-Sh:97}
\begin{eqnarray}
\mathbf{Var}\big\{\mathcal{N}_{n}[\varphi ]\big\} &=&\frac{1}{2}\int_{%
\mathbb{R}}\int_{\mathbb{R}}\left\vert \varphi (\lambda _{1})-\varphi
(\lambda _{2})\right\vert ^{2}K_{n}^{2}(\lambda _{1},\lambda _{2})d\lambda
_{1}d\lambda _{2}  \label{VK} \\
&=&\int_{\mathbb{R}}\int_{\mathbb{R}}\left\vert \frac{\Delta \varphi }{%
\Delta \lambda }\right\vert ^{2}\mathcal{V}_{n}(\lambda _{1},\lambda
_{2})d\lambda _{1}d\lambda _{2},  \label{Var_n}
\end{eqnarray}%
where%
\begin{equation}
K_{n}(\lambda _{1},\lambda _{2})=\sum_{l=0}^{n-1}\psi _{l}^{(n)}(\lambda
_{1})\psi _{l}^{(n)}(\lambda _{2})  \label{Kn}
\end{equation}%
is known as the reproducing kernel of the system (\ref{psi}),%
\begin{equation}
\frac{\Delta \varphi }{\Delta \lambda }=\frac{\varphi (\lambda _{1})-\varphi
(\lambda _{2})}{\lambda _{1}-\lambda _{2}},  \label{dfdl}
\end{equation}%
and
\begin{equation}
\mathcal{V}_{n}(\lambda _{1},\lambda _{2})=\Big[r_{n-1}^{(n)}\Big(\psi
_{n}^{(n)}(\lambda _{1})\psi _{n-1}^{(n)}(\lambda _{2})-\psi
_{n-1}^{(n)}(\lambda _{1})\psi _{n}^{(n)}(\lambda _{2})\Big)\Big]^{2}\Big/2.
\label{Vll}
\end{equation}%
Note that in passing from (\ref{VK}) to (\ref{Var_n}) we used the
Christoffel-Darboux formula \cite{Sz:75}%
\begin{equation}
K_{n}(\lambda _{1},\lambda _{2})=r_{n-1}^{(n)}\frac{\psi _{n}^{(n)}(\lambda
_{1})\psi _{n-1}^{(n)}(\lambda _{2})-\psi _{n-1}^{(n)}(\lambda _{1})\psi
_{n}^{(n)}(\lambda _{2})}{\lambda _{1}-\lambda _{2}}.  \label{CD}
\end{equation}%
It is of interest to have a spectral-theoretic interpretation of the above
formulas. Notice that $J^{(n)}$ of (\ref{J_n}) determines a selfadjoint
operator, that we denote again $J^{(n)}$. It acts in $l^{2}(\mathbb{Z}_{+})$%
, and the matrix
\begin{equation}
\mathcal{E}^{(n)}(d\lambda )=e^{(n)}(\lambda )d\lambda ,\;e^{(n)}(\lambda
)=\{e_{lm}^{(n)}(\lambda )\}_{l,m=0}^{\infty },\;e_{lm}^{(n)}(\lambda )=\psi
_{l}^{(n)}(\lambda )\psi _{m}^{(n)}(\lambda )  \label{Een}
\end{equation}%
is its resolution of identity \cite{Ak:62}. This allows us to write (\ref%
{Vll}) in the form%
\begin{eqnarray}
\mathcal{V}_{n}(\lambda _{1},\lambda _{2})
&=&(r_{n-1}^{(n)})^{2}(e_{n,n}^{(n)}(\lambda _{1})e_{n-1,n-1}^{(n)}(\lambda
_{2})  \label{Vee} \\
&+&e_{n,n}^{(n)}(\lambda _{2})e_{n-1,n-1}^{(n)}(\lambda
_{1})-2e_{n-1,n}^{(n)}(\lambda _{1})e_{n,n-1}^{(n)}(\lambda _{2}))/2.  \notag
\end{eqnarray}%
Assume now that $\varphi $ is of the class $C^{1}$ and%
\begin{equation}
\sup_{\lambda \in \mathbb{R}}\left\vert \varphi ^{\prime }(\lambda
)\right\vert <\infty .  \label{Ho}
\end{equation}%
for some positive $C$. It follows then from (\ref{orth}), (\ref{Var_n}) -- (%
\ref{Vll}) that%
\begin{equation}
\mathbf{Var}\big\{\mathcal{N}_{n}[\varphi ]\big\}\leq
(r_{n-1}^{(n)})^{2}\left( \sup_{\lambda \in \mathbb{R}}\left\vert \varphi
^{\prime }(\lambda )\right\vert \right) ^{2}.  \label{VCr}
\end{equation}%
It can be shown \cite{Pa-Sh:97} that if the potential satisfies the local
Lipshitz condition (\ref{Lipl}), then the coefficients $r_{n-1}^{(n)}$ are
bounded in $n$. We conclude that the variance of a linear statistic (\ref{LS}%
) is bounded in $n$ if the test function satisfies (\ref{Ho}):%
\begin{equation}
\mathbf{Var}\big\{\mathcal{N}_{n}[\varphi ]\big\}\leq \mathrm{Const}\left(
\sup_{\lambda \in \mathbb{R}}\left\vert \varphi ^{\prime }(\lambda
)\right\vert \right) ^{2}.  \label{Vbou}
\end{equation}%
This has to be compared with the well known fact of probability theory,
according to which the variance of a linear statistics of i.i.d. random
variables is $O(n)$ for any test function.

Notice that the condition (\ref{Ho}) for (\ref{Vbou}) to be valid is
pertinent. Indeed, it is shown in Appendix A.1 that for the Gaussian Unitary
Ensemble (\ref{GUEg}) and for the piece-wise constant $\varphi _{n}=\chi
_{\Delta }$ in (\ref{LS}), i.e., for (\ref{ECM}), the variance is $O(\log n)$%
, as $n\rightarrow \infty $ (see formulas (\ref{NDND}), (\ref{Ctou}) -- (\ref%
{Cins})).


\subsection{Asymptotic behavior of $\protect\psi _{l}^{(n)}(\protect\lambda)$%
}

\noindent We present now powerful asymptotic formulas for orthonormal
functions $\psi _{l}^{(n)}$ of (\ref{psi}) -- (\ref{orth}) due to Deift et
al \cite{De-Co:99} , valid in the case of a real analytic potential $V$ in (%
\ref{PdM}). We give a bit different their form that will be more convenient
below. Our form is reminiscent of standard semi-classical formulas and is
similar to that, given in \cite{Bl-It:99,Da-Co:00}. Notice that \cite%
{Da-Co:00} contains another, although heuristic, derivation of the
asymptotic formulas as well as an interesting heuristic explanation of the
quasiperiodicity of formulas for $q\geq 2$, based on ideas of statistical
mechanics and quantum field theory.

Assume that $V$ is real analytic. Then the support (\ref{sig}) of $N$ is a
union of $q<\infty $ finite disjoint intervals \cite{De-Co:99}. The function
$P$ in (\ref{rho}), (\ref{P}) is also real analytic. Following \cite%
{De-Co:99} we say that $V$ is regular if the inequality (\ref{Vns}) is
strict and $P$ is strictly positive on $\sigma $. Hence the DOS $\rho $ of (%
\ref{DOS}) is strictly positive inside $\sigma $ and vanishes precisely like
a square root at the endpoints. Denote%
\begin{equation}
N(\lambda )=N([\lambda ,\infty )).  \label{NNl}
\end{equation}%
According to \cite{De-Co:99}, in the regular case there exist functions $%
d_{n}(\lambda )$, and $\gamma _{n}(\lambda )$ such that if $\lambda $
belongs to the interior of the support (\ref{sig}), then
\begin{equation}
\psi _{n}^{(n)}(\lambda )=(2d_{n}(\lambda ))^{1/2}\cos \big(\pi nN(\lambda
)+\gamma _{n}(\lambda )\big)+O\left( n^{-1}\right) ,\ n\rightarrow \infty .
\label{as}
\end{equation}%
Moreover, $d_{n}(\lambda )$ and $\gamma _{n}(\lambda )$ depend on $n$ via
the vector $n\beta $, where $\beta =(\beta _{1},..,\beta _{q-1})$ is given
by (\ref{bl}). This means that there exist $n$-independent continuous
functions $\mathcal{D}:\sigma \times \mathbb{T}^{q-1}\rightarrow \mathbb{R}%
_{+}$, and $\mathcal{G}:\sigma \times \mathbb{T}^{q-1}\rightarrow \mathbb{R}$
such that
\begin{equation}
d_{n}(\lambda )=\mathcal{D}(\lambda ,n\beta ),\;\gamma (\lambda )=\mathcal{G}%
(\lambda ,n\beta ).  \label{dDgG}
\end{equation}%
If $\lambda $ belongs to the exterior of $\sigma $, then $\psi
_{n}^{(n)}(\lambda )$ decays exponentially in $n$ as $n\rightarrow \infty $.

Similar asymptotic formulas are valid for coefficients of the Jacobi matrix $%
J^{(n)}$ of (\ref{J_n}), i.e., there exist $n$-independent continuous
functions $\mathcal{R}:\mathbb{T}^{q-1}\rightarrow \mathbb{R}_{+}$ and $%
\mathcal{S}:\mathbb{T}^{q-1}\rightarrow \mathbb{R}$ such that
\begin{equation}
r_{n-1}^{(n)}=\mathcal{R}(n\beta )+O\left( n^{-1}\right) ,\;s_{n}^{(n)}=%
\mathcal{S}(n\beta )+O\left( n^{-1}\right) ,\;n\rightarrow \infty .
\label{rR}
\end{equation}%
The functions $\mathcal{D},\mathcal{G},\mathcal{R}$, and $\mathcal{S}$ can
be expressed via the Riemann theta-function, associated in the standard way
with two-sheeted Riemann surface obtained by gluing together two copies of
the complex plane slit along the gaps $%
(b_{1},a_{2}),...,(b_{q-1},a_{q}),(b_{q},a_{1})$ of the support of the
measure $N$, the last gap goes through the infinity.

The components of the vector $\beta =\{\beta _{l}\}_{l=1}^{q-1}$ are
rationally independent generically in $V$, thus the functions $\mathcal{D}%
(\lambda ,n\beta ),\mathcal{G}(\lambda ,n\beta ),\mathcal{R}(n\beta )$, and $%
\mathcal{S}(n\beta )$ are quasiperiodic in $n$ in general.

We will also need asymptotic formulas for $\psi _{n+k}^{(n)}$ as $%
n\rightarrow \infty $ and $k=O(1)$ (in particular, we need the case $k=-1$
in (\ref{Vll})). They can be extracted from \cite{De-Co:99} (see also \cite%
{Da-Co:00})), but it will be convenient to present a different derivation,
because the corresponding argument will be used in our analysis of limiting
laws for linear eigenvalue statistics. To this end we replace a regular
potential $V$ in (\ref{p_n}) by $V/g,\;g>0$, introducing explicitly the
amplitude of the potential. Then the quantities of asymptotic formulas (\ref%
{NNl}) -- (\ref{rR}) will depend on $g$, and it follows from the results of
\cite{De-Co:99,Ku-Mc:00} that these quantities will be continuous functions
of $g$ in a certain neighborhood of $g=1$, provided that the support (\ref%
{sig}) for $g=1$ consists of $q$ disjoint intervals. Consider now $%
r_{n+k-1}^{(n)}(g)$. Taking into account that the origin of the super-index $%
n$ in the above formulas is the factor $n$ in front of $V$ in (\ref{wn}), we
can write%
\begin{equation}
n\frac{V}{g}=(n+k)\frac{V}{g(1+k/n)}.  \label{Vng}
\end{equation}%
In other words, to obtain $r_{n+k-1}^{(n)}(g)$ and $\psi
_{n+k}^{(n)}(\lambda ,g)$ we have to make the change
\begin{equation}
g\rightarrow g+\frac{k}{n}g  \label{gren}
\end{equation}%
in the inverse amplitude of the potential. We obtain in view of (\ref{rR})
for $k=o(n)$:
\begin{eqnarray}
r_{n+k-1}^{(n)}(g) &=&r_{n+k-1}^{(n+k)}((1+k/n)g)  \label{rnR} \\
&\simeq &\mathcal{R}\left( (1+k/n)g,(n+k)\beta \left( (1+k/n)g\right) \right)
\notag \\
&\simeq &\mathcal{R}\left( g,n\beta (g)+k\alpha (g)\right) ,  \notag
\end{eqnarray}%
where%
\begin{equation}
\alpha (g)=(g\beta (g))^{\prime },  \label{agb}
\end{equation}%
and the symbol "$\simeq $" denotes here and below the leading term(s) of the
corresponding l.h.s. as $n\rightarrow \infty $.

We have an analogous formula for (\ref{NNl}):
\begin{equation*}
(n+k)N\left( \lambda ,(1+k/n)g\right) \simeq nN(\lambda ,g)+k\nu (\lambda
,g),
\end{equation*}%
where%
\begin{equation}
\nu (\lambda ,g)=\frac{\partial }{\partial g}(gN(\lambda ,g)).  \label{nudN}
\end{equation}%
By using these formulas, we can write for any fixed $k$ (in fact $k=o(n)$):
\begin{eqnarray}
\psi _{n+k}^{(n)}(\lambda ,g) &\simeq &(2\mathcal{D}(\lambda ,g,n\beta
+k\alpha ))^{1/2}  \label{pnnk} \\
&\times &\cos \Big(\pi nN(\lambda ,g)+\pi k\nu (\lambda ,g)+\mathcal{G}%
(\lambda ,g,n\beta +k\alpha )\Big).  \notag
\end{eqnarray}%
Applying to (\ref{rnR}) and its analog for $s_{n+k}^{(n)}$ the limiting
procedure (\ref{bnl}), we obtain the coefficients
\begin{equation}
r_{k-1}(x)=\mathcal{R(}x+k\alpha \mathcal{)},\;s_{k}(x)=\mathcal{S(}%
x+k\alpha \mathcal{)},\;x\in \mathbb{T}^{q-1},\;k\in \mathbb{Z}.
\label{rksk}
\end{equation}%
They determine a family of the double infinite Jacobi matrices $J(x),\;x\in
\mathbb{T}^{q-1}$:
\begin{equation}
(J(x)\psi )_{k}=r_{k}(x)\psi _{k+1}+s_{k}(x)\psi _{k}+r_{k-1}(x)\psi
_{k-1},\ \ k\in \mathbb{Z},  \label{J}
\end{equation}%
that can be viewed as a quasiperiodic operator, acting in $l^{2}(\mathbb{Z})$%
.

Consider now the functional%
\begin{equation}
\mathcal{E}[m]=-\int_{\sigma }\log |\lambda -\mu |m(d\lambda )m(d\mu ),
\label{Epot}
\end{equation}%
defined on unit non-negative measures, whose support is contained in $\sigma
$. This is the standard variational problem of the potential theory \cite%
{Sa-To:97}. Denote $\nu $ the unique minimizer of (\ref{Epot}). Then,
according to \cite{Bu-Ra:99} (see also \cite{Bu-Pa:02,Pa:96}), the
non-increasing function $\nu (\lambda )=\nu ([\lambda ,\infty ))$ (cf (\ref%
{NNl}))\ coincides with the function, defined in (\ref{nudN}). Moreover,
according to \cite{Pa:05} the measure $\nu $ is the Integrated Density of
States (IDS) measure of $J(x)$ (see \cite{Pa-Fi:92} for a definition of the
IDS measure in a general setting of ergodic operators), the support (\ref%
{sig}) is spectrum of $J(x)$, and (cf (\ref{bl}))%
\begin{equation}
\alpha _{l}=\nu ([a_{l},\infty )),\;l=1,...,q-1  \label{anu}
\end{equation}%
are the frequencies of quasi-periodic coefficients (\ref{rksk}) of $J(x)$.
In other words, $J(x)$ is a "finite band" Jacobi matrix, well known in
spectral theory and integrable systems \cite{Te:00}.

By using these facts and also taking into account that $\psi
_{n}^{(n)}(\lambda )$ decays exponentially in $n$ outside the support \cite%
{De-Co:99}, we can prove that for any continuous $\Phi :\mathbb{R}%
\rightarrow \mathbb{C}$ of a compact support we have in the limit (\ref{bnl}%
) of (\ref{Een})
\begin{equation}
\lim_{n_{j}(x)\rightarrow \infty }\int_{\mathbb{R}}\Phi (\lambda
)e_{n_{j}(x)+l,n_{j}(x)+m}^{(n)}(\lambda )d\lambda =\int_{\mathbb{\sigma }%
}\Phi (\lambda )e_{lm}(\lambda ,x)d\lambda ,  \label{eF}
\end{equation}%
where we have for $\lambda \in \sigma $:%
\begin{equation}
e_{lm}(\lambda ,x)=\psi _{l}^{+}(\lambda ,x)\overline{\psi _{m}^{+}(\lambda
,x)}+\psi _{l}^{-}(\lambda ,x)\overline{\psi _{m}^{-}(\lambda ,x)},
\label{epsi}
\end{equation}
\begin{equation}
\psi _{l}^{+}(\lambda ,x)=e^{i\pi l\nu (\lambda )}\mathcal{U}(\lambda
,x+l\alpha ),\;\mathcal{U}(\lambda ,x)=(\mathcal{D}(\lambda ,x)/2)^{1/2}e^{i%
\mathcal{G}(\lambda ,x)},  \label{psiD}
\end{equation}%
and $\psi _{l}^{-}=\overline{\psi _{l}^{+}\text{ }}$. The above formula can
also be written as%
\begin{equation}
\psi _{l}^{\pm }(\lambda ,x)=\psi _{0}^{\pm }(\lambda ,x+T^{l}\alpha ),
\label{psi0}
\end{equation}%
where $T:\mathbb{T}^{q-1}\rightarrow \mathbb{T}^{q-1}$ is defined as $%
Tx=x+\alpha $. This shows that $\Psi ^{\pm }(\lambda ,x)=\{\psi _{l}^{\pm
}(\lambda ,x)\}_{l\in \mathbb{Z}}$ is a generalized (quasi-Bloch)
eigenfunction of $J(x)$:%
\begin{equation*}
J(x)\Psi ^{\pm }(\lambda ,x)=\lambda \Psi ^{\pm }(\lambda ,x).
\end{equation*}%
In fact, $\{\Psi ^{\pm }(\lambda ,x)\}_{\lambda \in \sigma }$ form a
complete system in $l^{2}(\mathbb{Z})$).

We refer the reader to \cite{Pa:05} for more details of this aspect of
asymptotic formulas of \cite{De-Co:99}.

\subsection{Asymptotic behavior of variance}

\noindent We assume in this subsection that the test function $\varphi $ in (%
\ref{LS}) is of the class $C^{1}$ and does not depend on $n$. Hence,
function (\ref{dfdl}) is continuous in $(\lambda _{1},\lambda _{2})$. As a
result, fast oscillating in $n$ functions, entering (\ref{pnnk}), do not
contribute to the limit, as was already in obtaining (\ref{eF}) -- (\ref%
{psiD}). We have then from (\ref{Var_n}), (\ref{Vee}), and (\ref{eF}) -- (%
\ref{psiD})
\begin{equation}
\mathbf{Var}\big\{\mathcal{N}_{n}[\varphi ]\big\}\simeq \mathcal{V}(n\beta ),
\label{Vas}
\end{equation}%
where
\begin{equation}
\mathcal{V}(x)=\int_{\mathbb{\sigma }}\int_{\mathbb{\sigma }}\left\vert
\frac{\Delta \varphi }{\Delta \lambda }\right\vert ^{2}\mathcal{V}(\lambda
_{1},\lambda _{2},x)d\lambda _{1}d\lambda _{2},  \label{Vx}
\end{equation}
\begin{equation}
\mathcal{V}(\lambda _{1},\lambda _{2},x)=\mathcal{R}^{2}(x)\left(
e_{0,0}(\lambda _{1},x)e_{-1,-1}(\lambda _{2},x)-e_{-1,0}(\lambda
_{1},x)e_{-1,0}(\lambda _{2},x)\right) ,  \label{Vxl}
\end{equation}%
and $e_{lm}(\lambda ,x)$ are given by (\ref{epsi}) -- (\ref{psiD}), in
particular%
\begin{equation}
\int_{\sigma }e_{lm}(\lambda ,x)d\lambda =\delta _{lm}.  \label{en}
\end{equation}%
Since the charges $\beta _{l},\;l=1,...,q-1$ of (\ref{bl}) are continuous
and non-constant functionals of the potential (see \cite{De-Co:99,Ku-Mc:00}%
), the leading term of variance is a quasi-periodic function generically in
the potential. In particular, it has no limit as $n\rightarrow \infty $. Its
limiting points are indexed by the subset $\mathbb{H}^{q-1}\subset \mathbb{T}%
^{q-1}$, the closure of limiting points of the $q-1$ dimensional vectors (%
\ref{bn}). $\mathbb{H}^{q-1}$ is $\mathbb{T}^{q-1}$ generically in $V$, but
it can also be a proper subset of $\mathbb{T}^{q-1}$ (see, e.g. (\ref{bBP})).

The simplest case is $q=1$ of a single interval support. Here $\mathbb{H}%
^{0} $ is a point, there exist the limits \cite{APS:97,De-Co:98}%
\begin{equation}
\lim_{n\rightarrow \infty }r_{n+k-1}^{(n)}=r,\;\lim_{n\rightarrow \infty
}s_{n+k}^{(n)}=s,\;\forall k\in \mathbb{Z},  \label{1intr}
\end{equation}%
and the "limiting" Jacobi matrix $J(x)$ of (\ref{J}) has constant
coefficients:
\begin{equation*}
J_{jk}=r\delta _{j+1,k}+s\delta _{j,k}+r\delta _{j-1,k}.
\end{equation*}%
Placing the origin of the spectral axis at $s$, we obtain that $\lambda
=2r\cos \pi \nu $, $\psi ^{\pm }$ are just the plane waves,
\begin{equation*}
\sigma =[-2r,2r],
\end{equation*}%
and%
\begin{equation}
\mathcal{D}(\lambda )=-\nu ^{\prime }(\lambda )=\frac{1}{\pi \sqrt{%
4r^{2}-\lambda ^{2}}},\text{ }\lambda \in \sigma \text{.}  \label{Dnu}
\end{equation}%
This and general formula (\ref{Vxl}) yield a version of (\ref{Vas}) -- (\ref%
{Vx}) in which the role of $\mathcal{V}(x,\lambda _{1},\lambda _{2})$ plays
\begin{equation}
\mathcal{V}^{(1)}(\lambda _{1},\lambda _{2})=\frac{1}{4\pi ^{2}}\frac{%
4r^{2}-\lambda _{1}\lambda _{2}}{\sqrt{4r^{2}-\lambda _{1}^{2}}\sqrt{%
4r^{2}-\lambda _{2}^{2}}},\;\lambda _{1},\lambda _{2}\in \sigma .
\label{1intV}
\end{equation}%
This form of the variance was first found in physical papers \cite%
{Am-Co:90,Br-Ze:93} and proved rigorously in \cite{Jo:98}. We see that in
the single interval case the variance is universal, i.e., its functional
form does not depend explicitly on the potential, the information on the
potential being encoded in the unique parameter $r$ of (\ref{1intr}). In
particular, we have (\ref{1intV}) for the Gaussian Unitary Ensemble (\ref%
{GUEg}) \cite{Ca:01}.

In the case (\ref{VBP}) -- (\ref{bBP}) $\mathbb{H}^{q-1}$ consists of $q$
points (see (\ref{bBP}) and the variance is a $q$-periodic function of $n$.
Example (\ref{vq2}) -- (\ref{rq2}) corresponds to the simplest non-trivial
case $q=2$, where $\beta _{1}=1/2,\;\mathbb{H}^{1}=\{0,1/2\}$, the matrix $%
J(x)$ is 2-periodic, its coefficients are
\begin{equation}
r_{k}=\frac{b-(-1)^{k}a}{2},\;s_{k}=0,  \label{rssym}
\end{equation}%
and the variance is asymptotically 2-periodic function in $n$ \cite%
{Da-Co:00,APS:01}:%
\begin{equation*}
\mathbf{Var}\big\{\mathcal{N}_{n}[\varphi ]\big\}\simeq \mathcal{V}%
^{(2)}(n/2),
\end{equation*}%
where $\mathcal{V}^{(2)}(x)$ is given by (\ref{Vx}) in which $\mathcal{V}%
(\lambda _{1},\lambda _{2},x),\;x\in \mathbb{H}^{1}=\{0,1/2\}$ is

\begin{eqnarray}
\mathcal{V}^{(2)}(\lambda _{1},\lambda _{2},x) &=&\frac{1}{2\pi ^{2}}\frac{%
\varepsilon _{\lambda _{1}}\varepsilon _{\lambda _{2}}}{\sqrt{|R_{2}(\lambda
_{1})|}\sqrt{|R_{2}(\lambda _{2})|}}  \label{2intV} \\
&\times &\left( (a^{2}-\lambda _{1}\lambda _{2})(b^{2}-\lambda _{1}\lambda
_{2})-(-1)^{2x}ab(\lambda _{1}-\lambda _{2})^{2},\right) ,\;\lambda
_{1},\lambda _{2}\in \sigma ,  \notag
\end{eqnarray}%
with
\begin{equation}
R_{2}(\lambda )=(\lambda ^{2}-a^{2})(\lambda ^{2}-b^{2}),  \label{R2sym}
\end{equation}%
and $\varepsilon _{\lambda }=1$ if $\lambda \in (-b,-a)$ and $\varepsilon
_{\lambda }=-1$ if $\lambda \in (a,b)$. In fact, these formulas are valid
for any real analytic and even potential, producing a symmetric two-interval
support
\begin{equation}
\sigma =[-b,-a]\bigcup \;[a,b],\;0<a<b<\infty ,  \label{2sym}
\end{equation}%
(see \cite{APS:01,Da-Co:00}). The general case of a two-interval and not
necessarily symmetric support was analyzed in \cite{Da-Co:00}, where it was
found that the variance can be expressed via the classical elliptic
functions of Jacobi and Weierstrass.

We conclude that a minimum modification of the limiting law of linear
eigenvalue statistics in the case of a multi-interval support of the IDS,
comparing with the case of i.i.d. random variables, could be a family of
normal laws, indexed by the points of $\mathbb{H}^{q-1}$. We shall see below
that this modification is not sufficient in certain cases.

We remark in conclusion of this section that formulas (\ref{Vas}) -- (\ref%
{Vxl}) allow us to characterize the universality classes of ensembles (\ref%
{PdM}) with respect to the variance in the global regime, i.e., the sets of
ensembles (potentials), leading to the same asymptotic form of the variance
of linear statistics in the regime. Namely, since the potential is present
in (\ref{Vas}) -- (\ref{Vxl}) only via the endpoints $(a_{1},...,b_{q})$ of
support and via the charges $(\beta _{1},...,\beta _{q-1})$ of all but one
intervals of the support, these parameters determine a universality class.
Notice that the parameters are not necessarily independent.


\section{Limiting laws}

\subsection{Laplace transform of the probability law of linear eigenvalue
statistics}

\noindent In this subsection we obtain an expression for the Laplace
transform of the probability law of linear eigenvalue statistics (\ref{LS})
via orthogonal polynomials. We consider here real-valued test functions $%
\varphi :\mathbb{R}\rightarrow \mathbb{R}$. The Laplace transform is
evidently%
\begin{equation}
Z_{n}[\varphi ]=\mathbf{E}_{V}\left\{ e^{-\overset{\circ }{\mathcal{N}}%
_{n}[\varphi ]}\right\} ,  \label{chf}
\end{equation}%
where $\mathbf{E}_{V}\left\{ ...\right\} $ denotes the expectation with
respect to (\ref{PdM}) (or (\ref{p_n})), determined by a given potential $V$%
, and
\begin{equation*}
\overset{\circ }{\mathcal{N}}_{n}[\varphi ]=\mathcal{N}_{n}[\varphi ]-%
\mathbf{E}_{V}\left\{ \mathcal{N}_{n}[\varphi ]\right\} .
\end{equation*}%
It is convenient to introduce the parameter$\;s\in \lbrack 0,1]$ and to
consider the function%
\begin{equation*}
F_{n}(s)=\log Z_{n}[s\varphi ],\;s\in \lbrack 0,1].
\end{equation*}%
It is easy to see that%
\begin{equation*}
F_{n}(0)=0,\;F_{n}^{\prime }(0)=-\mathbf{E}_{V}\left\{ \overset{\circ }{%
\mathcal{N}}_{n}[\varphi ]\right\} =0,
\end{equation*}%
and
\begin{equation*}
F_{n}^{\prime \prime }(s)=\mathbf{E}_{V+\frac{s\varphi }{n}}\left\{ \mathcal{%
N}_{n}^{2}[\varphi ]\right\} -\mathbf{E}_{V+\frac{s\varphi }{n}}^{2}\left\{
\mathcal{N}_{n}[\varphi ]\right\} :=\mathbf{Var}_{V+\frac{s\varphi }{n}%
}\left\{ \mathcal{N}_{n}[\varphi ]\right\} .
\end{equation*}%
This yields the following expression for the logarithm of (\ref{chf}):%
\begin{equation}
\log Z_{n}[\varphi ]=F_{n}(1):=F_{n}[\varphi ]=\int_{0}^{1}(1-s)\mathbf{Var}%
_{V+\frac{s\varphi }{n}}\left\{ \mathcal{N}_{n}[\varphi ]\right\} ds.
\label{logcf}
\end{equation}%
We mention that there exist another expression for the Laplace transform (%
\ref{chf}). It dates back to the Heine formulas in the theory of orthogonal
polynomials (see e.g. \cite{Sz:75}, Theorem 2.1.1) and can be easily
obtained from the Gram theorem:%
\begin{eqnarray*}
Z_{n}[\varphi ] &=&\det \left\{ \int_{\mathbb{R}}e^{-\overset{\circ }{%
\varphi }(\lambda )}\psi _{j}^{(n)}(\lambda )\psi _{k}^{(n)}(\lambda
)d\lambda \right\} _{j,k=1}^{n} \\
&=&e^{\mathbf{E}\{\mathcal{N}_{n}[\varphi] \}}\det \left( 1-K_{n,\varphi
}\right) ,
\end{eqnarray*}%
where%
\begin{equation*}
\overset{\circ }{\varphi }(\lambda )=\varphi (\lambda )-\mathbf{E}_{V}\{%
\mathcal{N}_{n}[\varphi ]\},\ \mathbf{E}_{V}\{\mathcal{N}_{n}[\varphi
]\}=n\int_{\mathbb{R}}\varphi (\lambda )\mathbf{E}_{V}\{N_{n}(d\lambda )\},
\end{equation*}%
and $K_{n,\varphi }$ is the integral operator, defined as%
\begin{equation*}
(K_{n,\varphi }f)(\lambda )=\int_{\mathbb{R}}K_{n}(\lambda ,\mu
)(1-e^{-\varphi (\mu )})f(\mu )d\mu ,\;\lambda \in \mathbb{R}.
\end{equation*}%
These formulas and their analogs for unitary matrices were used to prove
various versions of the Central Limit Theorem (see e.g. \cite%
{Ba:97,Jo:98,Ke-Na:00,So:00,So:02,Sp:87,Wi:02}).


\subsection{Asymptotic behavior of the Laplace transform}

\noindent We will assume in this subsection that $\varphi $ is real
analytic. According to (\ref{logcf}), (\ref{Var_n}), and (\ref{Vll}), we
have to find the asymptotic form of $\psi _{n}^{(n)}$ and $\psi _{n-1}^{(n)}$
for the potential $V+s\varphi /n$. We have already seen in the previous
section that adding terms of the order $O(n^{-1})$ to the potential we
obtain non-trivial contributions to the asymptotic formulas because of fast
oscillating in $n$ functions in the r.h.s. of (\ref{as}) -- (\ref{rR}), etc.
The $O(n^{-1})$ terms appeared there because of the passage $n\rightarrow
n+k $, leading to (\ref{agb}) -- (\ref{rksk}). In this case the terms are
proportional to the potential, since we just change its amplitude: $%
V\rightarrow V(1-k/n)$ (see (\ref{Vng}) -- (\ref{gren})). This required
derivatives (\ref{agb}) and (\ref{nudN}) of "frequencies" $\beta
_{l},\;l=1,...,q-1$, and $N(\lambda )$ of fast oscillating functions in (\ref%
{as}) -- (\ref{rR}) with respect to the inverse amplitude $g$ of the
potential.

On the other hand, to find the asymptotic behavior of the Laplace transform,
we have to add to the potential the term $s\varphi /n$ (see (\ref{logcf})).
Since $\varphi \neq V$ in general, this requires variational derivatives of
frequencies with respect to potential, i.e., we have to add the term $%
\varepsilon \varphi $ to the potential, and find the derivative of $\beta
_{l},\;l=1,...,q-1$, and $N(\lambda )$ with respect to $\varepsilon $ at $%
\varepsilon =0$.

Consider first the case $q=1$, where the support of the IDS is a single
interval. Here the dependence on $x$ of functions $\mathcal{D},\mathcal{G},%
\mathcal{R}$ and $\mathcal{S}$ of (\ref{as}) -- (\ref{rR}) is absent (see (%
\ref{1intr}) -- (\ref{1intV})). Hence the term $s\varphi /n$ is negligible
in the limit $n\rightarrow \infty $, because there are no fast oscillating
in $n$ functions in the asymptotics of $\psi
_{n+k}^{(n)},\;k=0,-1,\;r_{n-1}^{(n)}$ and $s_{n}^{(n)}$, and we obtain from
(\ref{logcf}) and (\ref{1intV}):%
\begin{equation}
\lim_{n\rightarrow \infty }F_{n}[\varphi ]=\lim_{n\rightarrow \infty }%
\mathbf{Var}\{\mathcal{N}_{n}[\varphi ]\}/2.  \label{clt1i}
\end{equation}%
Notice also that we have here the "genuine" limit as $n\rightarrow \infty $,
but not a sublimit (\ref{bnl}) along a subsequence. We conclude that the
Central Limit Theorem is valid in this case. This was proved in \cite{Jo:98}
by the variational method and for a rather broad class of potentials and
test functions (not necessarily real analytic).

As is was shown in the previous section, the variance of a linear statistics
with a $C^{1}$ test function has no limit as $n\rightarrow \infty $ if $%
q\geq 2$. Its sublimits are indexed by points of the "hull" $\mathbb{H}%
^{q-1}\subset \mathbb{T}^{q-1}$. Hence we cannot expect the traditional CLT (%
\ref{clt1i}), as in the case of $q=1$. Rather this should to be a collection
of the CLT, indexed by $\mathbb{H}^{q-1}$:%
\begin{equation}
\lim_{n_{j}(x)\rightarrow \infty }F_{n_{j}(x)}[\varphi
]=\lim_{n_{j}(x)\rightarrow \infty }\mathbf{Var}\{\mathcal{N}%
_{n_{j}(x)}[\varphi ]\}/2=\mathcal{V}(x)/2,\;x\in \mathbb{H}^{q-1},
\label{geCLT}
\end{equation}%
where $\{n_{j}(x)\}$ and $\mathcal{V}(x)$ are defined in (\ref{bnl}) and (%
\ref{Vx}) -- (\ref{Vxl}). We will call this the generalized CLT.

We will show now that the generalized CLT is not always the case for $q\geq
2 $. Recall that $N(\lambda )$ and $\beta _{l}$ are functionals of $V$ and
denote%
\begin{equation}
\overset{\cdot }{\beta }_{l}[\varphi ]=\left. \frac{\partial }{\partial
\varepsilon }\beta _{l}\right\vert _{\varepsilon =0},\;l=1,...,q-1,\;\;%
\overset{\cdot }{N}[\varphi ]=\left. \frac{\partial }{\partial \varepsilon }%
N(\lambda )\right\vert _{\varepsilon =0}  \label{bNdot}
\end{equation}%
the variational derivatives of $\beta _{l}$ and $N(\lambda )$ with respect
to $V$. $\overset{\cdot }{\beta }_{l}[\varphi ]$ and$\overset{\cdot }{\text{
}N}[\varphi ]$ are linear functionals of $\varphi $ and nonlinear
functionals of $V$. It follows from \cite{Ku-Mc:00} that they are well
defined if $V$ is real analytic and regular and $\varphi $ is real analytic
and such that $\max_{\lambda \in \mathbb{R}}|V(\lambda )/\varphi (\lambda
)|<\infty $.

Arguing as in Section 2.2, we obtain that in this case $\psi _{n+k}^{(n)}$
is given by (\ref{pnnk}) -- (\ref{rksk}) with the replacement
\begin{equation*}
k\alpha _{l}\rightarrow k\alpha _{l}+s\overset{\cdot }{\beta }_{l}[\varphi
],\;\pi k\nu \rightarrow \pi k\nu +\pi s\overset{\cdot }{N}[\varphi ].
\end{equation*}%
Now, assuming (\ref{bnl}) and taking into account (\ref{logcf}), (\ref{Var_n}%
), and (\ref{Vee}), we obtain from (\ref{logcf}):%
\begin{equation}
F[\varphi ]:=\lim_{n_{j}(x)\rightarrow \infty }\log Z_{n_{j}(x)}[\varphi
]=\int_{0}^{1}(1-s)\mathcal{V}(x+s\overset{\cdot }{\beta }[\varphi ])ds,
\label{limcf}
\end{equation}%
where $\mathcal{V}$ is given by (\ref{Vx}). According to (\ref{Vx}) $%
\mathcal{V}$ is a quadratic functional of $\varphi $. Hence the functional $%
F[\varphi ]$ is not quadratic in general, because of the presence of the
term $s\overset{\cdot }{\beta }[\varphi ]$ in the argument of the integrand
in (\ref{limcf}). In other words we have here a limiting law in the sense of
(\ref{bnl}), but the law is not necessarily Gaussian.

It seems that a general classification of possible cases is rather complex.
We thus will give several examples, showing different cases of asymptotic
behavior of the Laplace transform of the probability law of linear
eigenvalue statistics.

Consider first the case, where the test function is a multiple of the
potential:
\begin{equation}
\varphi (\lambda )=tV(\lambda ),\;t\in \mathbb{R}.  \label{fiV}
\end{equation}%
Then (\ref{agb}) and the relation $\overset{\cdot }{\beta }_{l}[V]=-\left.
\beta ^{\prime }(g)\right\vert _{g=1}$ yield:%
\begin{equation*}
\overset{\cdot }{\beta }_{l}[\varphi ]=\overset{\cdot }{\beta }%
_{l}[tV]=-t\left( \beta _{l}(1)-\alpha _{l}(1)\right) ,\;l=1,...,q-1,
\end{equation*}%
where $\alpha _{l},\;l=1,...,q-1$ are defined in (\ref{anu}). Hence, if
\begin{equation}
\beta _{l}(1)=\alpha _{l}(1),\;l=1,...,q-1,  \label{ab}
\end{equation}%
then the integrand in (\ref{limcf}) does not depend on $s$, and we obtain in
view of (\ref{Vas}) -- (\ref{Vx}) the generalized Central Limit Theorem (\ref%
{geCLT}).

The equality (\ref{ab}) is valid for any potential of the form (\ref{VBP})
with $g=1$ and $v^{2}-4$ having only simple and real zeros, because
according to \cite{Bu-Pa:02}
\begin{equation*}
\beta _{l}(1)=\alpha _{l}(1)=\frac{q-l}{q},\;l=1,...,q-1.
\end{equation*}%
It is also valid for any even potential, having two equal local minima and
one local maximum, which is high enough to produce a two-interval support (%
\ref{2sym}). In this case (\ref{ab}) results from the symmetry, implying
\begin{equation}
\beta _{1}=\alpha _{1}=1/2  \label{ba2}
\end{equation}%
(recall that in this case the vectors $\beta $ of (\ref{bl}) and $\alpha $
of (\ref{agb}) are one-dimensional: $\beta _{1}=N(a),\;\alpha _{1}=\nu (a)$).

In all these cases the limiting Jacobi matrix $J$ of (\ref{J}) is $q$%
-periodic ($q=2$ in the case of (\ref{2sym}), see (\ref{rssym})).

It can also be shown that we have the generalized Central Limit Theorem for
potentials (\ref{VBP}) and $\varphi =tv$ (here the limiting matrix $J$ is
also $q$-periodic).

To demonstrate a possibility to have a non-Gaussian limiting law, we
consider a simplest non-trivial case of even potential with the two-interval
support (\ref{2sym}) and of test function%
\begin{equation}
\varphi (\lambda )=t\lambda ,\;t\in \mathbb{R},  \label{linl}
\end{equation}%
i.e., the case of "linear" linear statistic%
\begin{equation}
t\sum_{l=1}^{n}\lambda _{l}^{(n)}=t\mathrm{Tr\;}M_{n}.  \label{LSlin}
\end{equation}%
Since in this case $\Delta \varphi /\Delta \lambda $ of (\ref{dfdl}) is
equal to $t$, it follows from (\ref{Vxl}) and (\ref{en}) that
\begin{equation}
\mathcal{V}(x)=t^{2}\mathcal{R}^{2}(x),  \label{VR}
\end{equation}%
and then (\ref{limcf}) implies that in the case (\ref{linl}) (and for any
support) we have for the exponent of the limiting Laplace transform:%
\begin{equation}
F[\varphi ]\Big|_{\varphi (\lambda )=t\lambda }=\lim_{n_{j}(x)\rightarrow
\infty }\log Z_{n_{j}(x)}[\varphi ]\Big|_{\varphi (\lambda )=t\lambda
}=t^{2}\int_{0}^{1}(1-s)\mathcal{R}^{2}\left( x+s\left. \overset{\cdot }{%
\beta }[\varphi ]\right\vert _{\varphi (\lambda )=t\lambda }\right) ds.
\label{Flin}
\end{equation}%
According to \cite{Da-Co:00}), it is possible to express the coefficient $%
\mathcal{R}(x)$, corresponding to the two-interval support, via the Jacobi
elliptic function:
\begin{equation}
\mathcal{R}^{2}(x)=\frac{(b-a)^{2}}{4}+\frac{ab}{2}\;\mathrm{cn}^{2}(x+1/2),
\label{Rcn}
\end{equation}%
where $\mathrm{cn}(x)=\mathrm{cn}(2K(k)x|k),\;k^{2}=4ab/(a+b)^{2},\;K(k)$ is
the elliptic integral of the first kind. In view of (\ref{ba2}) the
coefficient $r_{k}$ of (\ref{rksk}) is given by (\ref{rssym}):%
\begin{equation*}
r_{k-1}=\mathcal{R}\left( \frac{k}{2}\right) =\frac{b-(-1)^{k}a}{2},
\end{equation*}%
and is 2-periodic (see also \cite{APS:97}). In view of (\ref{VR}) this
implies that the variance of (\ref{LSlin}) is asymptotically 2-periodic in $%
n $:%
\begin{equation}
\mathcal{R}^{2}\left( \frac{n}{2}\right) =\frac{b^{2}+a^{2}}{4}-(-1)^{n}%
\frac{ab}{2}.  \label{valV}
\end{equation}%
Furthermore, it is shown in Appendix that
\begin{equation}
\left. \overset{\cdot }{\beta }_{1}[\varphi ]\right\vert _{\varphi (\lambda
)=t\lambda }=t\omega ,\;\omega =\frac{a}{4K(a/b)},  \label{bd2}
\end{equation}%
Hence we obtain from (\ref{Rcn})
\begin{equation}
\left. F[\varphi ]\right\vert _{\varphi (\lambda )=t\lambda
}=\int_{0}^{t}(t-s)\mathcal{R}^{2}(x+s\omega )ds.  \label{limq2}
\end{equation}%
It follows from (\ref{bd2}) that $\omega $ is irrational generically in $a$,
and $b$, hence $\mathcal{R}^{2}(x+s\omega )$ is quasi periodic in $s$ in
these cases. Since $\mathcal{R}^{2}$ is 1-periodic and real analytic, we can
write its Fourier series%
\begin{equation}
\mathcal{R}^{2}(x)=\sum_{m\in \mathbb{Z}}c_{m}e^{2\pi imx}  \label{RF}
\end{equation}%
with fast decaying coefficients. Plugging (\ref{bd2}), and (\ref{RF}) in (%
\ref{Flin}), we obtain%
\begin{equation}
\left. F[\varphi ]\right\vert _{\varphi (\lambda )=t\lambda }=\frac{%
c_{0}t^{2}}{2}-tA^{\prime }(x)-A(x)+A(x+\omega t),  \label{nGau}
\end{equation}%
where%
\begin{equation}
A(x)=\sum_{m\in \mathbb{Z}\backslash \{0\}}\frac{c_{m}}{(2\pi im\omega )^{2}}%
e^{2\pi imx}.  \label{Ax}
\end{equation}%
We see that the logarithm of the limiting Laplace transform of the
probability law of statistics (\ref{LSlin}) contains not only a multiple of $%
t^{2}/2$, that would correspond to the CLT, but also a linear in $t$ term, a
constant in $t$ term, and either quasi periodic (generically in $a,b$, when $%
\omega $ is irrational) or periodic (in special cases, where $\omega $ is
rational) function of $t$. Besides, while the variance of statistics (\ref%
{LSlin}) is (\ref{valV}) in the limit (\ref{bnl}), the coefficient in front
of $t^{2}/2$ is%
\begin{equation*}
c_{0}=\int_{\mathbb{T}}\mathcal{R}^{2}(x)dx=\frac{b^{2}+a^{2}}{4},
\end{equation*}%
hence is not the variance (\ref{valV}).

Notice also that in the case $q\geq 3$ we would have in (\ref{Ax}) the sum
over $\mathbb{Z}^{q-1}\backslash \{0\}$, and the expression%
\begin{equation*}
\overset{\cdot }{\beta }\cdot m:=\overset{\cdot }{\beta }_{1}m_{1}+...%
\overset{\cdot }{\beta }_{q-1}m_{q-1},
\end{equation*}%
in the denominator, that can be arbitrary small for certain collections of $%
(m_{1},...,m_{q-1})\in $ $\mathbb{Z}^{q-1}$ and $\overset{\cdot }{\beta }:=(%
\overset{\cdot }{\beta }_{1},...,\overset{\cdot }{\beta }_{q-1})$ $\in $ $%
\mathbb{R}^{q-1}$ in the case, where the components of $\overset{\cdot }{%
\beta }$ are irrational. Hence, to make the series in (\ref{Ax}) convergent,
we have to assume that the components of $\overset{\cdot }{\beta }$ are
sufficiently bad approximated by rationals (e.g. a Diophantine condition).

According to Appendix, in a general case of a real analytic $\varphi $ and a
two-interval support%
\begin{equation}
\sigma =[a_{1},b_{1}]\bigcup \;[a_{2},b_{2}],\;-\infty
<a_{1}<b_{1}<a_{2}<b_{2}<\infty ,  \label{2gen}
\end{equation}%
we have:%
\begin{equation}
\overset{\cdot }{\beta _{1}}[\varphi ]=-\frac{1}{2\pi iI}\int_{\sigma }\frac{%
\varphi (\mu )-\varphi (-\mu )}{\sqrt{R_{2}(\mu )}}d\mu ,  \label{bfI}
\end{equation}%
where $\sqrt{R_{2}(\mu )}$ is defined in (\ref{Xl}),
\begin{equation}
I=\int_{b_{1}}^{a_{2}}\frac{d\mu }{\sqrt{(b_{2}-\mu )(a_{2}-\mu )(\mu
-b_{1})(\mu -a_{1})}}=\frac{2}{((b_{2}-b_{1})(a_{2}-a_{1}))^{1/2}}K(\kappa ),
\label{I}
\end{equation}%
and%
\begin{equation*}
\kappa ^{2}=\frac{(a_{2}-b_{1})(b_{2}-a_{1})}{(b_{2}-b_{1})(a_{2}-a_{1})}
\end{equation*}%
(see \cite{Ry-Gr}, formula (3.149.4)). In the symmetric case (\ref{2sym}) we
have%
\begin{equation*}
I=\frac{2}{b+a}K\left( \frac{2\sqrt{ab}}{b+a}\right) =\frac{2}{a}K\left(
\frac{b}{a}\right) ,
\end{equation*}%
where the second equality results from the formula $(1+k)^{-1}K(2\sqrt{k}%
/(1+k))=K(k)$ \cite{Ry-Gr}, formula (8.126.3). It follows then that in this
case and for $\varphi =t\lambda $ (\ref{bfI}) coincides with (\ref{bd2}).

Combining (\ref{limcf}) for $q=2$ and (\ref{bfI}), we obtain the limiting
form of the Laplace transform of the probability law of linear eigenvalue
statistics for real analytic test functions. In fact, this form is also
valid for larger classes of test functions, in particular, for bounded $C^{1}
$ functions with bounded derivative. Indeed, for any such function $\varphi $
there exists a sequence $\{\varphi _{k}\}$ of real analytic test functions,
converging to $\varphi $ in the $C^{1}$ metrics. Besides:

(i) $\mathcal{V}$ of (\ref{Vx}) is continuous in $\varphi $ in the $C^{1}$
metrics and in $x\in \mathbb{T}^{1}$;

(ii) $\overset{\cdot }{\beta }_{1}$ of (\ref{bfI}) is continuous in $\varphi
$ in the $C$ metrics;

(iii) we \ have for the functional $F_{n}$ of (\ref{logcf})%
\begin{equation}
|F_{n}[\varphi _{1}]-F_{n}[\varphi _{2}]|\leq C\left( \sup_{\lambda \in
\mathbb{R}}|\varphi _{1}^{\prime }-\varphi _{2}^{\prime }|\right) ^{2},
\label{F2fb}
\end{equation}

where $C$ depend only on the potential and is finite if (\ref{Vlog}) holds.
To prove this inequality, implying the uniform in $n$ the$\ C^{1}$
continuity of $F_{n}$ in $\varphi $, we use an argument, similar to that
proving (\ref{logcf}), to obtain%
\begin{equation*}
F[\varphi _{1}]-F[\varphi _{2}]=\int_{0}^{1}(1-s)Var_{V+\varphi _{s}/n}\{%
\mathcal{N}_{n}[\varphi _{1}-\varphi _{2}]\}ds,
\end{equation*}%
where $\varphi _{s}=s\varphi _{1}+(1-s)\varphi _{2},\;s\in \lbrack 0,1]$.
This and (\ref{Vbou}) yield (\ref{F2fb}).

Now we can write:%
\begin{equation*}
|F_{n_{j}(x)}[\varphi ]-F_{x}[\varphi ]|\leq |F_{n_{j}(x)}[\varphi
]-F_{n_{j}(x)}[\varphi _{k}]|+|F_{n_{j}(x)}[\varphi ]-F_{x}[\varphi
_{k}]|+|F_{x}[\varphi ]-F_{x}[\varphi _{k}]|
\end{equation*}%
for $F_{n}$ of (\ref{logcf})\ and $F_{x}$ of (\ref{limcf}) with explicitly
indicated dependence on $x\in \mathbb{T}^{1}$. Making the limit $%
n_{j}(x)\rightarrow \infty $ and using (\ref{F2fb}), we obtain in view of
the above result on real analytic test functions:%
\begin{equation*}
\lim_{n_{j}(x)\rightarrow \infty }\sup |F_{n_{j}(x)}[\varphi ]-F_{x}[\varphi
]|\leq C\left( \sup_{\lambda \in \mathbb{R}}|\varphi ^{\prime }-\varphi
_{k}^{\prime }|\right) ^{2}+|F_{x}[\varphi ]-F_{x}[\varphi _{k}]|.
\end{equation*}%
Now the $C^{1}$ limit $\varphi _{k}\rightarrow \varphi $ and above
assertions (i) and (ii) yield the validity of our results (\ref{limcf}) for $%
q=2$ and (\ref{bfI}) for bounded $C^{1}$ test functions with bounded
derivative, hence a possibility to have non-Gaussian limiting laws (the
non-validity of the Central Limit Theorem) for linear eigenvalue statistics
of ensembles (\ref{PdM}) with these test functions.

Formula (\ref{bfI}) allows us to characterize the class of potentials and
test functions for which the (generalized) Central Limit Theorem (\ref{geCLT}%
) is valid in the case of a general two-interval support (\ref{2gen}).
Indeed, for any pair $(V,\varphi )$ for which the r.h.s. of (\ref{bfI}) is
zero the integrand of (\ref{limcf}) does not depend on $s$ and we have the
generalized CLT (\ref{geCLT}). In particular, in the symmetric case (\ref%
{2sym}) it follows from (\ref{bfI}) that $\overset{\cdot }{\beta _{1}}%
[\varphi ]$ is zero if and only if
\begin{equation*}
\int_{a}^{b}\frac{\varphi (\mu )-\varphi (-\mu )}{\sqrt{(b^{2}-\mu ^{2})(\mu
^{2}-a^{2})}}d\mu =0,
\end{equation*}%
In particular, for an even potential of support (\ref{2sym}) and an even
test function $\varphi $ the generalized CLT is valid.

One can view $\varphi $ as an analog of external field in statistical
mechanics. Hence, we can say that in this case an even external field "does
not break the symmetry". On the other hand a "generic" $\varphi $ or an odd $%
\varphi $, such that
\begin{equation*}
\int_{a}^{b}\frac{\varphi (\mu )d\mu }{\sqrt{(b^{2}-\mu ^{2})(\mu ^{2}-a^{2})%
}}\neq 0,
\end{equation*}%
is a "breaking symmetry field" and leads to a non-Gaussian limiting law. Its
simplest case $\varphi (\lambda )=t\lambda $ (\ref{linl}) is given by (\ref%
{nGau}).

\section{Intermediate and local regimes}

\noindent In this section we consider limiting laws of linear eigenvalue
statistics for the test functions, given by (\ref{rint}) \ and (\ref{rloc})
with a $C^{1}$ function $\varphi $ and $\lambda _{0}$ belonging to the
interior of $\sigma $.

We begin again by calculating the asymptotic form of the variance in these
cases. Changing variables to%
\begin{equation}
\lambda _{1,2}=\lambda _{0}+t_{1,2}/n^{\alpha },\;0<\alpha \leq 1,
\label{lt}
\end{equation}%
we obtain from (\ref{Var_n}):%
\begin{equation}
\mathbf{Var}\big\{\mathcal{N}_{n}[\varphi ]\big\}=\int_{\mathbb{R}}\int_{%
\mathbb{R}}\Big(\frac{\Delta \varphi }{\Delta t}\Big)^{2}\mathcal{V}%
_{n}(\lambda _{0}+t_{1}/n^{\alpha },\lambda _{0}+t_{2}/n^{\alpha
})dt_{1}dt_{2}.  \label{Vta}
\end{equation}%
To find the asymptotic form of the r.h.s. we will use again (\ref{Vee}) and (%
\ref{Een}) in which $\psi _{n+k}^{(n)}(\lambda ),\;k=0.-1$ is given by (\ref%
{pnnk}) with$\ \lambda =\lambda _{0}+t/n^{\alpha }$. Taking into account
that $\mathcal{D},N$, and $\mathcal{G}$ are smooth functions of $\lambda $
in a sufficiently small neighborhood of $\lambda _{0}$, we can write
\begin{eqnarray*}
\psi _{n+k}^{(n)}(\lambda ) &\simeq &\left( 2\mathcal{D}(\lambda _{0},n\beta
+k\alpha )\right) ^{1/2} \\
&\times &\cos \left( \pi nN(\lambda _{0})+\pi k\nu (\lambda _{0})-\pi \rho
(\lambda _{0})n^{1-\alpha }+\mathcal{G}(\lambda _{0},n\beta +k\nu )\right) ,
\end{eqnarray*}%
where $\rho (\lambda )=-N^{\prime }(\lambda )$, and we do not indicate the
dependence on $g$ (in fact it suffices to consider the case $g=1$). Plugging
this into (\ref{Een}) and (\ref{Vee}), and omitting in the resulting
integrand of the r.h.s. of (\ref{Vta}) the fast oscillating terms, we obtain
\begin{eqnarray}
\mathbf{Var}\big\{\mathcal{N}_{n}[\varphi ]\big\} &\simeq &B(\lambda
_{0},n\beta )  \label{Vni} \\
&\times &\int_{\mathbb{R}}\int_{\mathbb{R}}(\varphi (t_{1})-\varphi
(t_{2}))^{2}\frac{\sin ^{2}\left( \pi \rho (\lambda
_{0})(t_{1}-t_{2})n^{1-\alpha }\right) }{2\pi ^{2}(t_{1}-t_{2})^{2}}%
dt_{1}dt_{2},  \notag
\end{eqnarray}%
where%
\begin{eqnarray}
B(\lambda ,x) &=&2\pi ^{2}\mathcal{R}^{2}(x)\mathcal{D}(\lambda ,x)\mathcal{D%
}(\lambda ,x-\alpha )  \label{Bx} \\
&&\times \sin ^{2}\left( \pi \nu (\lambda )+\mathcal{G}(\lambda ,x)-\mathcal{%
G}(\lambda ,x-\alpha )\right) .  \notag
\end{eqnarray}%
This leads to the following result in the local regime $\alpha =1$ and for
the limit (\ref{bnl}):%
\begin{eqnarray}
\lim_{n_{j}(x)\rightarrow \infty }\mathbf{Var}\{\mathcal{N}%
_{n_{j}(x)}[\varphi _{n_{j}(x)}]\} &=&B(\lambda _{0},x)  \label{VlocB} \\
&\times &\int_{\mathbb{R}}\int_{\mathbb{R}}(\varphi (t_{1})-\varphi
(t_{2}))^{2}\frac{\sin ^{2}\left( \pi \rho (\lambda
_{0})(t_{1}-t_{2})\right) }{2\pi ^{2}(t_{1}-t_{2})^{2}}dt_{1}dt_{2}.  \notag
\end{eqnarray}%
It is known that in the local regime the variance of linear eigenvalue
statistics has a universal limiting form \cite{De-Co:99,Pa-Sh:97}
\begin{equation}
\lim_{n\rightarrow \infty }\mathbf{Var}\{N_{n}[\varphi _{n}]\}=\int_{\mathbb{%
R}}\int_{\mathbb{R}}(\varphi (t_{1})-\varphi (t_{2}))^{2}\frac{\sin
^{2}\left( \pi \rho (\lambda _{0})(t_{1}-t_{2})\right) }{2\pi
^{2}(t_{1}-t_{2})^{2}}dt_{1}dt_{2},  \label{Vloc}
\end{equation}%
in which all the information on the potential is encoded in $\rho (\lambda
_{0}).$ This and (\ref{VlocB}) imply the identity%
\begin{equation}
B(\lambda ,x)=1,\;  \label{id}
\end{equation}%
valid for all $x\in \mathbb{T}^{q-1}$ and $\lambda $, belonging to the
interior of the $\sigma $ (a direct proof of the identity can be extracted
from the proof of Lemma 6.1 of \cite{De-Co:99}).

In the intermediate regime $0<\alpha <1$ we still have a fast oscillating
factor
\begin{equation*}
\sin ^{2}\left( n^{1-\alpha }\pi \rho (\lambda _{0})(t_{1}-t_{2})\right)
\end{equation*}%
in the integrand of (\ref{Vni}). Replacing the factor by its average $1/2$,
and using (\ref{id}) we obtain in this regime%
\begin{equation}
\lim_{n\rightarrow \infty }\mathbf{Var}\{N_{n}[\varphi _{n}]\}=\int_{\mathbb{%
R}}\int_{\mathbb{R}}\frac{(\varphi (t_{1})-\varphi (t_{2}))^{2}}{4\pi
^{2}(t_{1}-t_{2})^{2}}dt_{1}dt_{2}.  \label{Vint}
\end{equation}%
As in the case of the local regime the limit here is the same for any
subsequence (\ref{bnl}), hence we do not need to assume (\ref{bnl}). We
conclude that the variance of linear statistics has a well defined limit in
the intermediate regime as well. Moreover, (\ref{Vint}) is the "smoothed"
version of variance (\ref{Vloc}) in the local regime, since (\ref{Vint}) is (%
\ref{Vloc}) in which the "oscillating" factor $\sin ^{2}(\pi \rho (\lambda
_{0})(t_{1}-t_{2}))$ is replaced by its average $1/2$.

Passing to the Fourier transform
\begin{equation*}
\widehat{\varphi }(k)=\frac{1}{2\pi }\int_{\mathbb{R}}e^{ikt}\varphi (t)dt,
\end{equation*}%
we can rewrite the r.h.s. of (\ref{Vint}) in the form
\begin{equation*}
\int_{\mathbb{R}}|k|\widehat{\varphi }(k)\widehat{\varphi }(-k)dk
\end{equation*}%
appearing in the continuous analog of the strong Szeg\"{o} theorem (see \cite%
{Ba:97,Di-Ev:01,So:00} and references therein).

The universality property of unitary invariant Matrix Models \cite%
{De-Co:99,Pa-Sh:97} implies that the Laplace transform of the probability
law of linear eigenvalue statistics has the following limiting form in the
local regime%
\begin{equation}
\left. \lim_{n\rightarrow \infty }Z_{n}[\varphi _{n}]\right\vert _{\varphi
_{n}(\lambda )=\varphi ((\lambda -\lambda _{0})n)}=e^{2\pi \rho (\lambda
_{0})\widehat{\varphi }(0)}\det (1-S_{\varphi }),  \label{Zloc}
\end{equation}%
where $S_{\varphi }$ is the integral operator, defined as
\begin{equation}
(S_{\varphi }f)(t)=\int_{\mathbb{R}}\frac{\sin \pi \rho (\lambda _{0})(t-u)}{%
\pi \rho (\lambda _{0})(t-u)}(1-e^{-\varphi (u)})f(u)du,  \label{Sphi}
\end{equation}%
and we assume that $\varphi $ in (\ref{rloc}) is continuous and integrable
on $\mathbb{R}$. It is obvious that the logarithm of the r.h.s. of (\ref%
{Zloc}) is not quadratic in $\varphi $, hence the CLT is not valid in the
local regime (see e.g. \cite{Hu-Ru:03} for related results).

If, however, we take in the above formulas
\begin{equation*}
\varphi (t)=\Phi ((t-t_{0})\delta ),
\end{equation*}%
where $\Phi $ does not depend on $\delta \rightarrow 0$, and $t_{0}\in
\mathbb{R}$, i.e., we assume that the test function in (\ref{Zloc}) -- (\ref%
{Sphi}) is "slow varying", then it can be shown (see e.g. \cite%
{Ba:97,Co-Le:95,Sp:87}) that the limit of the r.h.s. of (\ref{Zloc}) as $%
\delta \rightarrow 0$ is the r.h.s. of (\ref{Vint}), divided by 2.

On the other hand, take as a test function in (\ref{limcf})%
\begin{equation}
\varphi (\lambda )=\Phi ((\lambda -\lambda _{0})/\delta ),  \label{fqint}
\end{equation}%
where $\Phi $ does not depend on $\delta \rightarrow 0$, and $\lambda _{0}$
belongs to the interior of the support of $N$, i.e., assume that $\varphi $
is "fast varying". Since the variational derivatives (linear functionals of $%
\varphi $) $\overset{\cdot }{\beta }_{l}[\varphi ],\;l=1,...,q-1$ of (\ref%
{bNdot}) can be written as%
\begin{equation*}
\overset{\cdot }{\beta }_{l}[\varphi ]=\int_{\mathbb{R}}b_{l}(\lambda
)\varphi (\lambda )d\lambda ,
\end{equation*}%
we have
\begin{equation*}
\overset{\cdot }{\beta }_{l}[\varphi ]=\delta \int_{\mathbb{R}}b_{l}(\lambda
_{0}+\delta t)\Phi (t)dt\rightarrow 0,\;\delta \rightarrow 0.
\end{equation*}%
Hence the term $s\overset{\cdot }{\beta }$ in the argument of the integrand
of (\ref{limcf}) vanishes in the limit $\delta \rightarrow 0$ and we obtain
from (\ref{Vx}), changing variables to $\lambda _{1,2}=\lambda _{0}+\delta
t_{1,2}$:%
\begin{equation}
\left. \lim_{\delta \rightarrow \infty }F[\varphi ]\right\vert _{\varphi
(\lambda )=\Phi ((\lambda -\lambda _{0})/\delta )}=\frac{1}{2}\mathcal{V}(x),
\label{limeF}
\end{equation}%
where%
\begin{equation*}
\mathcal{V}(x)=\mathcal{V}(\lambda _{0},\lambda _{0},x)\int_{\mathbb{R}%
}\int_{\mathbb{R}}\frac{(\Phi (t_{1})-\Phi (t_{2}))^{2}}{2\pi
^{2}(t_{1}-t_{2})^{2}}dt_{1}dt_{2}.
\end{equation*}%
Now it can be shown, by using (\ref{Vxl}), (\ref{psiD}), (\ref{psi0}), and (%
\ref{Bx}), that the r.h.s. of the last formula coincides with (\ref{Vint}).
Hence the limit (\ref{limeF}) coincides again with the r.h.s. of (\ref{Vint}%
), divided by 2.

The above suggests that the CLT is valid in the intermediate regime. This
was indeed proved in several particular cases (see \cite{So:00} and
references therein).

\renewcommand{\thesection}{\Alph{section}} \setcounter{section}{1} %
\setcounter{equation}{0}

\appendix

\section{Appendices}

\subsection{Variance of the Eigenvalue Counting Measure of the GUE}

It has been shown in Sections 2.1 and 2.3 that the variance of linear
eigenvalue statistics is of the order $O(1)$ if the corresponding test
function is of the class $C^{1}$. Here we will argue that this condition is
close to be optimal. To this end we will show that linear eigenvalue
statistics for the GUE (\ref{GUEg}) are of the order $O(\log n)$ if $\varphi
$ is the indicator $\chi _{\Delta }$ of an interval $\Delta =(a,b)\in
(-2w,2w)$. For this choice of $\varphi $ the corresponding linear statistic
is the Eigenvalue Counting Measure (\ref{ECM}). We will derive the
asymptotic formula
\begin{equation}
\mathbf{Var}\{\mathcal{N}_{n}(\Delta )\}=\frac{1}{\pi ^{2}}\log
n\;(1+o(1)),\ n\rightarrow \infty .  \label{NDND}
\end{equation}%
We have from (\ref{VK})

\begin{equation}
\mathbf{Var}\{\mathcal{N}_{n}(\Delta )\}=\int_{a}^{b}d\lambda \int_{\mathbb{R%
}\backslash (a,b)}K_{n}^{2}(\lambda ,\mu )d\mu ,  \label{VKD}
\end{equation}%
where $K_{n}$ is given by (\ref{Kn}). According to (\ref{CD}) we need
asymptotics of $\psi _{l}^{(n)}$ for $l=n-1,n$ to find the asymptotic
behavior of (\ref{VKD}). These follow from (\ref{pnnk}) (with the remainder $%
o(1)$), (\ref{Dnu}), and (\ref{GUEg}) or can be obtained from the well known
Plancherel-Rotah formulas for the Hermite polynomials (see \cite{Sz:75},
formula (8.22.12)):%
\begin{equation}
\psi _{n+k}^{(n)}(\lambda )=\frac{1}{(2\pi \sqrt{g}\sin \theta )^{1/2}}\cos
\left( n\alpha (\theta )+k\theta +\mathcal{G}(\theta )\right) (1+o(1)),
\label{PRved}
\end{equation}%
where $\lambda =2\sqrt{g}\cos \theta ,\;\alpha (\theta )=\theta -\sin
2\theta /2,\;\mathcal{G}(\theta )=\alpha (\theta )/2-\pi /4.$

Besides, in this case \cite{Sz:75}%
\begin{equation}
r_{n-1}^{(n)}=\sqrt{g}.  \label{rln}
\end{equation}%
We write (\ref{VKD}) as
\begin{equation}
\mathbf{Var}\{\mathcal{N}_{n}(\Delta )\}=I_{1}+I_{2}+I_{3}+I_{4},  \label{I4}
\end{equation}%
where
\begin{eqnarray}
I_{1} &=&\int_{b-a}^{\infty }dx\int_{0}^{b-a}K_{n}^{2}(b-y,b-y+x)dy,\;
\label{I1I2} \\
I_{2} &=&\int_{0}^{b-a}dx\int_{0}^{x}K_{n}^{2}(b-y,b-y+x)dy,  \notag
\end{eqnarray}%
and $I_{3}$ and $I_{4}$ can be obtained from $I_{1}$ and $I_{2}$ by
replacing $K_{n}(b-y,b-y+x)$ by $K_{n}(a+y,a+y-x)$.

By using the Christoffel-Darboux formula (\ref{CD}) we can write the
integrands in $I_{1,2}$ as
\begin{equation*}
(r_{n-1}^{(n)})^{2}\Psi _{n}^{2}(x,y)x^{-2},
\end{equation*}%
where%
\begin{equation}
\Psi _{n}(x,y)=\left. \left( \psi _{n}^{(n)}(\lambda )\psi _{n-1}^{(n)}(\mu
)-\psi _{n-1}^{(n)}(\lambda )\psi _{n}^{(n)}(\mu )\right) \right\vert
_{\lambda =b-y,\mu =b-y+x}.  \label{Ppsi}
\end{equation}%
In view of (\ref{rln}) the integral $I_{1}$ is bounded from above by
\begin{equation*}
g(b-a)^{-2}\int_{\mathbb{R}^{2}}\left( \psi _{n}^{(n)}(\lambda )\psi
_{n-1}^{(n)}(\mu )-\psi _{n-1}^{(n)}(\lambda )\psi _{n}^{(n)}(\mu )\right)
^{2}d\lambda d\mu =2g/(b-a)^{2},
\end{equation*}%
where in writing the last equality we used the orthonormality of the system $%
\{\psi _{l}^{(n)}\}_{l\geq 0}$.

To find the asymptotic behavior of $I_{2}$ we will use formula (\ref{PRved}%
). However its direct use is impossible since the remainder term leads to
the divergent integral in $x$. Thus, we write this integral as the sum of
the two, over the intervals $(0,A/n)$ and $(A/n,b-a)$, where $A$ is fixed.
In the first integral we write (\ref{Ppsi}) as
\begin{equation*}
\Psi _{n}(x,y)=\left. \left( \psi _{n}^{(n)}(\lambda )-\psi _{n}^{(n)}(\mu
)\right) \psi _{n-1}^{(n)}(\mu )-\left( \psi _{n-1}^{(n)}(\lambda )-\psi
_{n-1}^{(n)}(\mu )\right) \psi _{n}^{(n)}(\mu )\right\vert _{\lambda =b-y,\
\mu =b-y+x}.
\end{equation*}%
In view of (\ref{PRved}) $\psi _{n+k}^{(n)},\;k=0,-1$ are of order $O(1)$.
Besides, we have the relation \cite{Sz:75}%
\begin{equation*}
\frac{d}{d\lambda }\psi _{l}^{(n)}=-\frac{n}{g}\lambda \psi
_{l}^{(n)}+\left( \frac{nl}{g}\right) ^{1/2}\psi _{l-1}^{(n)},
\end{equation*}%
implying that $\frac{d}{d\lambda }\psi _{n+k}^{(n)}$ is $O(n),\;k=0,-1$.
This yields the bound $\left\vert \Psi _{n}(x,y)\right\vert \leq \mathrm{%
const}\cdot nx$, according to which the first integral is%
\begin{equation*}
\mathrm{const}\cdot n^{2}\int_{0}^{A/n}xdx=O(1),\;n\rightarrow \infty .
\end{equation*}%
In the integral over $(A/n,b-a)$ we use the asymptotic formula (\ref{pnnk}),
neglecting, as it was done in the proof of the previous sections, the fast
oscillating terms, and noticing that the remainder term $o(1)$ in (\ref{pnnk}%
) yields the error $o(\log n)$. This and (\ref{rln}) for $l=n$\ lead to the
asymptotic expression (cf \ref{1intV})
\begin{equation*}
\frac{1}{2\pi ^{2}}\int_{A/n}^{b-a}\frac{dx}{x^{2}}\int_{0}^{x}dy\left.
\frac{4g-\lambda \mu }{\sqrt{4g-\lambda ^{2}}\sqrt{4g-\mu ^{2}}}\right\vert
_{\lambda =b-y,\ \mu =b-y+x}+o(1)
\end{equation*}%
The leading contribution to the integral is due the integral in $x$ over an
interval $(A/n,\varepsilon )$, where $\varepsilon >0$ is small enough and $n$%
-independent, because the integral over $(\varepsilon ,b-a)$ is $n$%
-independent, hence is $O(1),\;n\rightarrow \infty $. Then the condition $%
\lambda =b-y,\ \mu =b-y+x$ can be replaced by $\lambda =b-y,\ \mu =b-y$. As
a result the integral in $y$ is asymptotically $x$, and we obtain
\begin{equation*}
I_{1}+I_{2}=\frac{1}{2\pi ^{2}}\log n+O(1),\;n\rightarrow \infty .
\end{equation*}%
The same contribution is due the sum $I_{3}+I_{4}$ in (\ref{I4}), hence we
obtain (\ref{NDND}).

Analogous formula is also known for unitary matrices. For this and
corresponding Central Limit Theorem see \cite{Ke-Na:00,Wi:02}. \bigskip

Consider the covariance $\mathbf{Cov}\{\mathcal{N}_{n}(\Delta _{1}),\mathcal{%
N}_{n}(\Delta _{2})\}$ of the Eigenvalue Counting Measures (\ref{ECM}) for $%
\Delta _{1}=(a_{1},b_{1})$ and $\Delta _{2}=(a_{2},b_{2}),\;\Delta
_{1,2}\subset (-2\sqrt{g},2\sqrt{g})$. It follows from (\ref{p_n}) that (cf (%
\ref{VK}))%
\begin{eqnarray*}
&&\mathbf{Cov}\{\mathcal{N}_{n}[\varphi _{1}],\mathcal{N}_{n}[\varphi
_{2}]\}:=\mathbf{E}\{\mathcal{N}_{n}[\varphi _{1}]\mathcal{N}_{n}[\varphi
_{2}]\}-\mathbf{E}\{\mathcal{N}_{n}[\varphi _{1}]\}\mathbf{E}\{\mathcal{N}%
_{n}[\varphi _{2}]\} \\
&=&\frac{1}{2}\int_{\mathbb{R}}\int_{\mathbb{R}}\Big(\varphi _{1}(\lambda
_{1})-\varphi _{1}(\lambda _{2})\Big)\Big(\varphi _{2}(\lambda _{1})-\varphi
_{2}(\lambda _{2})\Big)K_{n}^{2}(\lambda _{1},\lambda _{2})d\lambda
_{1}d\lambda _{2}.
\end{eqnarray*}%
This yields for $\varphi _{1,2}=\chi _{\Delta _{1,2}}$
\begin{equation}
\mathbf{Cov}\{\mathcal{N}_{n}(\Delta _{1}),\mathcal{N}_{n}(\Delta
_{2})\}=\int_{\mathbb{R}^{2}}\left[ \chi _{(\Delta _{1}\cap \Delta
_{2})\times \mathbb{R}}(\lambda ,\mu )-\chi _{\Delta _{1}\times \Delta
_{2}}(\lambda ,\mu )\right] K_{n}^{2}(\lambda ,\mu )d\lambda d\mu .
\label{Covchi}
\end{equation}%
The argument, proving (\ref{NDND}), allows us to find the leading terms of
the covariance as $n\rightarrow \infty $ in various cases. We will assume
below for the sake of definiteness that $a_{1}<b_{1},\;a_{2}<b_{2},\;b_{1}%
\leq b_{2},\;$and $b_{1}-a_{1}\leq b_{2}-a_{2}$, and consider the following
cases of the asymptotic behavior of (\ref{Covchi}) (see \cite{Wi:02} for
analogous results for unitary matrices).

\noindent (i). Disjoint intervals ($b_{1}<a_{2}$):%
\begin{equation}
-\frac{1}{4\pi ^{2}}\int_{\Delta _{1}}d\lambda \int_{\Delta _{2}}\frac{%
4g-\lambda \mu }{\sqrt{4g-\lambda ^{2}}\sqrt{4g-\mu ^{2}}}d\mu ,
\label{Cdis}
\end{equation}%
i.e., in this case the covariance is bounded as it was for the case of
linear statistics, generated by $C^{1}$ functions in (\ref{Vbou}) and (\ref%
{Vas}).

\noindent (ii). Touching (outside) intervals ($a_{2}=b_{1}$):%
\begin{equation}
-\frac{1}{2\pi ^{2}}\log n.  \label{Ctou}
\end{equation}

\noindent (iii) Touching (inside) intervals ($a_{1}=a_{2}$, but $b_{1}<b_{2}$%
):%
\begin{equation}
\frac{1}{2\pi ^{2}}\log n.  \label{Cins}
\end{equation}

\noindent (iv) Embedded intervals ($a_{2}<a_{1},\;b_{1}<b_{2}$):%
\begin{equation}
\frac{1}{4\pi ^{2}}\int_{\Delta _{1}}d\lambda \int_{\mathbb{R}\backslash
\Delta _{2}}\frac{4g-\lambda \mu }{\sqrt{4g-\lambda ^{2}}\sqrt{4g-\mu ^{2}}}%
d\mu .  \label{Cemb}
\end{equation}

\noindent (v) Intersecting intervals ($a_{2}<b_{1}$ but $a_{1}<a_{2}$ and $%
b_{1}<b_{2}$). In this case we can write the equality
\begin{equation}
\chi _{(\Delta _{1}\cap \Delta _{2})\times \mathbb{R}}-\chi _{\Delta
_{1}\times \Delta _{2}}=\chi _{(a_{1},a_{2})\times (\mathbb{R}\setminus
(a_{1},b_{2}))}-\chi _{(a_{1},a_{2})\times (b_{1},b_{2})},  \label{2chi}
\end{equation}%
showing that the domain of integration in (\ref{Covchi}) does not contain
the \textquotedblleft dangerous\textquotedblright\ diagonal $\lambda =\mu $.
Hence, the leading term of the covariance in this case is of the order $O(1)$%
, and the respective coefficient is given by the integral of the product of
the integrand of (\ref{Cemb}) and of (\ref{2chi}).

\medskip By using similar reasoning, it is possible to find the asymptotic
form of the variance and the covariance of the Eigenvalue Counting Measure
in other cases. Note that in the "regular" cases, where the covariance is
bounded, (see (\ref{Vas}) and (\ref{Cdis})), its leading term is additive in
$\Delta _{1}$ and $\Delta _{2}$ (or in $\varphi _{1}$ and $\varphi _{2}$),
while in the "singular" cases, where $\mathbf{Cov}\{\mathcal{N}_{n}(\Delta
_{1}),\mathcal{N}_{n}(\Delta _{2})\}=O(\log n)$, its leading term is
independent of $\Delta _{1}$ and $\Delta _{2}$.

\subsection{Variational derivative of frequency in the two-interval case}

\noindent We derive here formula (\ref{bd2}). We will use the variational
approach, based on the functional (\ref{EV}).

Write the minimum condition (\ref{Vs}), (\ref{Vef}) for $V+\varepsilon
\varphi $, and compute its derivative at $\varepsilon =0$. This yields%
\begin{equation}
\varphi (\lambda )-2\int_{\sigma }\log |\lambda -\mu |\overset{\cdot }{\rho }%
(\mu )d\mu =\mathrm{const,\;}\lambda \in \sigma ,  \label{ELf}
\end{equation}%
where%
\begin{equation*}
\overset{\cdot }{\rho }=\left. \frac{\partial }{\partial \varepsilon }\rho
\right\vert _{\varepsilon =0},
\end{equation*}%
and $\rho $ is the density of the measure $N$ of (\ref{IDS}). Notice that
the differentiation of the limits of integration in (\ref{Vef}) does not
contribute to (\ref{ELf}), because $\rho $ vanishes at each endpoints of the
support according to (\ref{rho}).

The derivative of (\ref{ELf}) in $\lambda $ is the singular integral
equation (cf (\ref{inteq})):%
\begin{equation*}
\mathrm{v.p.}\int_{\sigma }\frac{\overset{\cdot }{\rho }(\mu )d\mu }{\mu
-\lambda }=-\frac{\varphi ^{\prime }(\lambda )}{2},\;\lambda \in \sigma .
\end{equation*}%
The general solution of the equation in the case (\ref{2sym}) is \cite{Mu}%
\begin{equation*}
\frac{C_{1}\lambda +C}{X_{2}(\lambda )}+\frac{1}{2\pi ^{2}X_{2}(\lambda )}\;%
\mathrm{v.p.}\int_{\sigma }\frac{\varphi ^{\prime }(\mu )X_{2}(\mu )d\mu }{%
\mu -\lambda },
\end{equation*}%
where%
\begin{equation*}
X_{2}(\lambda )=-i\sqrt{R_{2}(\lambda )},\;\sqrt{R_{2}(\lambda )}=\left.
\sqrt{R_{2}(z)}\right\vert _{z=\lambda +i0},
\end{equation*}%
$R_{2}(z)=(z^{2}-a^{2})(z^{2}-b^{2})$ \ (see (\ref{Xl}) for $q=2$) and $%
\sqrt{R_{2}(z)}$ is the branch of the square root, fixed by the condition: $%
\sqrt{R(z)}=z^{2}+O(z),\;z\rightarrow \infty $. The branch assumes pure
imaginary values of opposite sides on the edges of $\sigma $, seen as a cut
of $\mathbb{C}$.

Taking into account the equalities%
\begin{equation*}
\int_{\sigma }\overset{\cdot }{\rho }(\mu )d\mu =0,
\end{equation*}%
and%
\begin{equation*}
\int_{\sigma }\frac{d\mu }{\sqrt{R_{2}(\mu )}}=0,\;\int_{\sigma }\frac{\mu
d\mu }{\sqrt{R_{2}(\mu )}}=-\pi i,
\end{equation*}%
we find that
\begin{equation}
\overset{\cdot }{\rho }(\lambda )=\frac{iC}{\sqrt{R_{2}(\lambda )}}+\frac{1}{%
2\pi ^{2}\sqrt{R_{2}(\lambda )}}\mathrm{v.p.}\int_{\sigma }\frac{\varphi
^{\prime }(\mu )\sqrt{R_{2}(\mu )}d\mu }{\mu -\lambda }.  \label{rdfi}
\end{equation}%
The constant $C$ can be found as follows. Denote $f(z)$ the Stieltjes
transform of $\rho $:%
\begin{equation}
f(z)=\int_{\sigma }\frac{\rho (\mu )d\mu }{\mu -z},\;z\notin \sigma .
\label{f}
\end{equation}%
Recalling that $V$ is real analytic, using (\ref{rho}), (\ref{P}), and
arguing as above, we find that%
\begin{equation}
f(z)=-\frac{V^{\prime }(z)}{2}-\frac{\sqrt{R_{2}(z)}}{2\pi i}\int_{\sigma }%
\frac{V^{\prime }(\mu )-V^{\prime }(z)}{\mu -z}\frac{d\mu }{\sqrt{R_{2}(\mu )%
}}  \label{fV}
\end{equation}%
(the analog of the formula with $R_{q}$ instead of $R_{2}$ is valid for any
finite number $q$ of intervals in (\ref{sig})).

Write now the minimum condition (\ref{Vs}) as
\begin{equation*}
V_{eff}(b)-V_{eff}(a)=0,
\end{equation*}%
or, in view of (\ref{Vef}) and (\ref{f}), as%
\begin{equation}
\int_{-a}^{a}\left( f(\lambda )+\frac{V^{\prime }(\lambda )}{2}\right)
d\lambda =0.  \label{gapc}
\end{equation}%
The condition is valid for any potential, in particular, for $V+\varepsilon
\varphi $. According to (\ref{fV}) for any sufficiently small $\varepsilon $
the integrand here is proportional to $\sqrt{R_{2}(\lambda )}$, in which $a$
and $b$ are now functions of $\varepsilon $. Hence, the integrand vanishes
at the edges of the support and the derivative of (\ref{gapc}) with $V$
replaced by $V+\varepsilon \varphi $ with respect to $\varepsilon $ at $%
\varepsilon =0$ is%
\begin{equation}
\int_{-a}^{a}\left( \overset{\cdot }{f}(\lambda )+\frac{\varphi ^{\prime
}(\lambda )}{2}\right) d\lambda =0.  \label{gapd}
\end{equation}%
This and the formula%
\begin{equation}
\overset{\cdot }{f}(z)=\int_{\sigma }\frac{\overset{\cdot }{\rho }(\mu )d\mu
}{\mu -z}=-\frac{\pi C}{\sqrt{R_{2}(z)}}-\frac{1}{2\pi i\sqrt{R_{2}(z)}}%
\int_{\sigma }\frac{\varphi ^{\prime }(\mu )\sqrt{R_{2}(\mu )}}{\mu -z}d\mu ,
\label{fdfi}
\end{equation}%
following from (\ref{rdfi}), yield in the case $\varphi (\lambda )=t\lambda $%
:%
\begin{equation*}
\overset{\cdot }{\rho }(\lambda )=\frac{C}{\sqrt{(b^{2}-\lambda
^{2})(\lambda ^{2}-a^{2})}}+\frac{P_{2}(\lambda )}{2\pi \sqrt{(b^{2}-\lambda
^{2})(\lambda ^{2}-a^{2})}},\;\lambda \in \sigma ,
\end{equation*}%
where $P_{2}(\lambda )=\lambda ^{2}-(a^{2}+b^{2})/2$,
\begin{equation}
C=\frac{tI_{2}}{2\pi I_{1}},  \label{C}
\end{equation}%
and%
\begin{equation}
I_{1}=\int_{-a}^{a}\frac{d\lambda }{\sqrt{(b^{2}-\lambda ^{2})(a^{2}-\lambda
^{2})}},\;I_{2}=2\int_{-a}^{a}\frac{P_{2}(\lambda )d\lambda }{\sqrt{%
(b^{2}-\lambda ^{2})(a^{2}-\lambda ^{2})}}.  \label{I12}
\end{equation}%
Now we can find $\overset{\cdot }{\beta _{1}}[\varphi ]$ for $\varphi
(\lambda )=t\lambda $. We have by (\ref{C}) and (\ref{I12}):%
\begin{equation}
\overset{\cdot }{\beta _{1}}:=\int_{a}^{b}\overset{\cdot }{\rho }(\lambda
)d\lambda =\frac{t}{2\pi I_{1}}(I_{2}J_{1}-I_{1}J_{2}),  \label{bdl}
\end{equation}%
where%
\begin{equation*}
J_{1}=\int_{a}^{b}\frac{d\lambda }{\sqrt{(b^{2}-\lambda ^{2})(\lambda
^{2}-a^{2})}},\;J_{2}=\int_{a}^{b}\frac{P_{2}(\lambda )d\lambda }{\sqrt{%
(b^{2}-\lambda ^{2})(\lambda ^{2}-a^{2})}}.
\end{equation*}%
By using standard formulas (see e.g. \cite{Ry-Gr}, formulas (3.159)), we find%
\begin{equation}
I_{1}=2K(k)/a,\;I_{2}=2a\left( \frac{1-k^{2}}{2}K(k)-E(k)\right) ,
\label{Is}
\end{equation}%
\begin{equation}
J_{1}=K(k^{\prime })/a,\;J_{2}=a\left( E(k^{\prime })-\frac{1+k^{2}}{2}%
K(k^{\prime })\right) ,  \label{Js}
\end{equation}%
where $K(k)$ and $E(k)$ are the complete elliptic integrals of the first and
second kind, $k=b/a,\;k^{\prime }=\sqrt{1-k^{2}}$. These formulas and the
identity $EK^{\prime }+E^{\prime }K-KK^{\prime }=\pi /2$, where $K^{\prime
}=K(k^{\prime }),\;E^{\prime }=E(k^{\prime })$ (\cite{Ry-Gr}, formula
(8.122)) imply (\ref{bd2}).

A more involved version of the above argument leads to (\ref{bfI}). We note
first that to prove the formula for a real analytic $\varphi $ it suffices
to consider
\begin{equation*}
\varphi _{z_{0}}(\lambda )=\frac{1}{\lambda -z_{0}},\;z_{0}\notin \sigma .
\end{equation*}%
We have in this case (see below):%
\begin{equation}
\overset{\cdot }{\beta _{1}}[\varphi _{z_{0}}]=-\frac{1}{2I\sqrt{R_{2}(z_{0})%
}},  \label{bz0}
\end{equation}%
where $I$ is defined in (\ref{I}).

Assuming that this formula is valid and using the Cauchy theorem to write%
\begin{equation*}
\varphi (\lambda )=\frac{1}{2\pi i}\int_{C_{\sigma }}\frac{\varphi
(z_{0})dz_{0}}{\lambda -z_{0}}, \; \lambda \in \sigma,
\end{equation*}%
where $C_{\sigma }$ is the contour encircling $\sigma $ in the clockwise
direction, we obtain in view of the linearity of $\overset{\cdot }{\beta _{1}%
}[\varphi ]$ in $\varphi $:
\begin{eqnarray*}
\overset{\cdot }{\beta _{1}}[\varphi ] &=&\frac{1}{2\pi i}\int_{C_{\sigma }}%
\overset{\cdot }{\beta _{1}}[\varphi _{z_{0}}]\varphi (z_{0})dz_{0} \\
&=&-\frac{1}{2I}\frac{1}{2\pi i}\int_{C_{\sigma }}\frac{\varphi (z_{0})}{%
\sqrt{R_{2}(z_{0})}}dz_{0}.
\end{eqnarray*}%
Now the relation $\sqrt{R_{2}(\lambda -i0)}=-\sqrt{R_{2}(\lambda +i0)}$
yields (\ref{bfI}).

To prove (\ref{bz0}) we use the general formulas (\ref{rdfi}), (\ref{fdfi}),
and (\ref{gapd}) with
\begin{equation*}
\varphi _{z_{0}}^{\prime }(\lambda )=-\frac{1}{(\lambda -z_{0})^{2}}=-\frac{%
\partial }{\partial z_{0}}\frac{1}{\lambda -z_{0}}.
\end{equation*}%
Arguing as in the case $\varphi (\lambda )=t\lambda $ above, we obtain an
analog of (\ref{bdl}) whose denominator contains $I$ of (\ref{I}) instead of
$I_{1}$ of (\ref{I12}), and whose numerator is a bilinear combination of
integrals of $|R_{2}(\mu )|^{-1/2}$ over $[a_{2},b_{2}]$ and of derivatives
with respect to $z_{0}$ of the integrals over $[b_{1},a_{2}]$ and $%
[a_{2},b_{2}]$ of $[(\mu -z_{0})|R_{2}(\mu )|^{1/2}]^{-1}$ and $P_{2}(\mu
)[(\mu -z_{0})|R_{2}(\mu )|^{1/2}]^{-1}$, where $P_{2}$ is defined now by
the relation $\sqrt{R_{2}(z)}=P_{2}(z)+O(1/z),\;z\rightarrow \infty $. These
integrals can be expressed via the complete elliptic integrals of the first,
second, and third kinds. Furthermore, the complete elliptic integrals of the
third kind can be expressed via the incomplete elliptic integrals of the
first and the second kinds, whose arguments depend on $z_{0}$. This allows
us to obtain a formula for $\overset{\cdot }{\beta _{1}}[\varphi _{z_{0}}]$,
whose numerator is expressed via the complete elliptic integrals of the
first and the second kind and derivatives with respect to $z_{0}$ of the
incomplete elliptic integrals of the first and the second kind. The
derivatives are proportional to $(R_{2}(z_{0}))^{-1/2}$ (see \cite{Ry-Gr},
formulas (8.123)) This and a bit tedious algebra lead eventually to (\ref%
{bz0}).

Another derivation of (\ref{bz0}) is given in \cite{Da-Co:00} (see formula
(3.14) of the paper). The derivation is based on a two step procedure of
minimization of (\ref{EV}): the first step is the minimization over all unit
measures with a given charge $\beta _{1}\in (0,1)$ of the "band" $%
[a_{2},b_{2}]$ of the support, and the second is the minimization of this
minimum over $\beta _{1}$.


\begin{thebibliography}{99}
\bibitem{Ak:62} Akhiezer, N., \ \emph{Classical Moment Problem} (Haffner,
New York,1965).

\bibitem{APS:97} Albeverio, S., Pastur, L., and Shcherbina, M., "On
asymptotic properties of certain orthogonal polynomials", \ Mathematical
Physics, Analysis, Geometry \ \textbf{4}, 263-277 (1997).

\bibitem{APS:01} Albeverio, S., Pastur, L., and Shcherbina, M., "On the $1/n$%
- expansion for some unitary invariant ensembles of random matrices", \
Commun. Math. Phys. \ \textbf{222}, 271-305 (2001).

\bibitem{Am-Co:90} Ambj\"{o}rn, J., Jurkiewicz, J., and Makeenko, Yu., \
"Multiloop correlators for two-dimensional quantum gravity", \ Phys. Lett. \
\textbf{B251}, 517-524 (1990).

\bibitem{Ba-Si:04} Bai, Z. D., Silverstein, J. W., "CLT for linear spectral
statistics of large dimensional sample covariance matrices", \ Ann. Prob. \
\textbf{32}, 553-605 (2004).

\bibitem{Ba:97} Basor, E., "Distribution functions for random variables for
ensembles of positive Hermitian matrices", \ Comm. Math. Phys. \ \textbf{188}%
, 327-350 (1997).


\bibitem{Bl-It:99} Bleher, P., Its, A., "Semiclassical asymptotics of
orthogonal polynomilas, Riemann-Hilbert problem, and universality of the
matrix models", \ Ann. Math. \ \textbf{150,} 185-266 (1999).

\bibitem{Da-Co:00} Bonnet, G., David. F., and Eynard, B., "Breakdown of
universality in multi-cut matrix models", \ J. Phys. \ \textbf{A33,}
6739-6768 (2000).

\bibitem{BPS:95} Boutet de Monvel, A., Pastur, L., and Shcherbina, M., "On
the statistical mechanics approach to the random matrix theory: the
integrated density of states", \ J. Stat. Phys. \textbf{79}, 585-611 (1995).

\bibitem{Br-Ze:93} Br\'{e}zin, E., and Zee, A., "Universality of the
correlations between eigenvalues of large random matrices", \ Nuclear Phys.
\textbf{B402 }, 613-627 (1993).

\bibitem{Bu-Pa:02} Buslaev, V., and Pastur, L., "A class of multi-interval
eigenvalue distributions of matrix models and related structures", \ in
\emph{Asymptotic Combinatorics and Applications to Mathematical Physics},
edited by A. Vershik, and V. Malyshev \ (Dordrecht, \ Kluwer, 2002) pp.
52-70.

\bibitem{Bu-Ra:99} Buyarov, V. S., and Rakhmanov, E. A., "Families of
equilibrium measures in an external field on the real axis", \ Sb. Math. \
\textbf{190}, 791-802 (1999).

\bibitem{Ca:01} Cabanal-Duvillard, T., "Fluctuations de la loi empirique de
grandes matrices al\'{e}atoires", \ Ann. Inst. H. Poincar\'{e}, Probab.
Statist. \ \textbf{37}, 373-402 (2001).

\bibitem{Ch-La:98} Chen, Y., Lawrence, N., "On the linear statistics of
Hermitian random matrices", J. Phys. \textbf{A34}, 1141-1152 (1998).

\bibitem{Co-Le:95} Costin. O., and Lebowitz, J. L., "Gaussian fluctuations
in random matrices", \ Phys. Rev. Lett. \ \textbf{75}, 69-72 (1995).

\bibitem{De-Co:98} Deift, P., Kriecherbauer, T., and McLaughlin, K., "New
results on the equilibrium measure for logarithmic potentials in the
presence of an external field", \ J. Approx. Theory \ \textbf{95}, 388-475
(1998).

\bibitem{De-Co:99} Deift, P., Kriecherbauer, T., McLaughlin, K., Venakides,
S., and Zhou, X., "Uniform asymptotics for polynomials orthogonal with
respect to varying exponential weights and applications to universality
questions in random matrix theory", \ Comm. Pure Appl. Math. \ \textbf{52},
1335-1425 (1999).

\bibitem{Di-Ev:01} Diaconis, P., and Evans, S., "Linear functionals of
eigenvalues of random matrices", \ Trans. of AMS \ \textbf{353}, 2615-2633
(2001).

\bibitem{Gi:01} Girko, V., \ \emph{Theory of Stochastic Canonical Equations}
\ (Dordrecht, Kluwer, 2001), vols. \ I, II.

\bibitem{Ry-Gr} Gradshteyn, I. S., and Ryzhik, I. M., \emph{Table of
Integrals, Series and Products} \ (San Diego, Academic Press, 1980).

\bibitem{Gu:05} Gustavsson, J., "Gaussian fluctuations of eigenvalues in the
GUE", Ann. Inst. H. Poincare: Prob. and Statistics, \textbf{41}, 151-178
(2005).

\bibitem{Hu-Ru:03} Hughes, C. P., Rudnick, Z., "Mock-Gaussian benavior for
linear statistics of classical compact groups", J. Phys. \textbf{A36},
2919-2932 (2003).

\bibitem{Jo:98} Johansson, K., "On fluctuations of eigenvalues of random
Hermitian matrices", \ Duke Math.J. \ \textbf{91}, 151-204 (1998).

\bibitem{Ke-Na:00} Keating, J. P., and Snaith, N., "Random matrix theory and
$\zeta (1/2+it)$", \ Commun. Math. Phys. \ \textbf{214}, 57-89 (2001).

\bibitem{KKP} Khoruzhenko, B., Khorunzhy, A., and Pastur, L., "$\ 1/n$%
-corrections to the Green functions of random matrices with independent
entries", \ J. of Phys. \textbf{A 28}, L31-L35 (1995).

\bibitem{Ku-Mc:00} Kuijlaars, A. B. J., and McLaughlin, K., "Generic
behaviour of the density of states in random matrix theory and equilibrium
problem of real analytic external field", \ Commun. Pure Appl. Math. \
\textbf{53}, 736-785 (2000).

\bibitem{Me:92} Mehta, M. L.,\ \emph{Random Matrices} \ (New York, Academic
Press, 1991).

\bibitem{Mu} Muskhelishvili, N. I., \ \emph{Singular Integral Equations} \
(Groningen, P. Noordhoff, 1953).

\bibitem{Pa:96} Pastur, L., "Spectral and probabilistic aspects of matrix
models", \ in \emph{Algebraic and Geometric Methods in Mathematical Physics}%
, \ edited by A. Boutet de Monvel, and V. Marchenko\ (Dodrecht, Kluwer,
1996) pp. 207-242.

\bibitem{Pa:00} Pastur, L., "Random matrices as paradigm", \ in \emph{%
Mathematical Physics 2000}, \ edited by A. Grigoryan, T. Kibble, and \ B.
Zegarlinski (London, Imperial College Press, 2000) pp. 216-266.

\bibitem{Pa:05} Pastur, L., "From random matrices to quasiperiodic Jacobi
matrices via orthogonal polynomials", \ J. Approx. Theory \ \textbf{139},
269-292 (2006).

\bibitem{Pa-Fi:92} Pastur, L., and Figotin, A., \ \emph{Spectra of Random
and Almost Periodic Operators} \ (Berlin, Springer, 1992).

\bibitem{Pa-Sh:97} Pastur, L., and Shcherbina, M., "Universality of the
local eigenvalue statistics for a class of unitary invariant matrix
ensembles", \ J. Stat. Phys. \ \textbf{8}, 109-147 (1997).

\bibitem{Sa-To:97} Saff, E. B., and Totik, V., \ \emph{Logarithmic
Potentials with External Fields} \ (New York, Springer, 1997).

\bibitem{Si-So:97} Sinai, Ya., and Soshnikov, A., "Central limit theorem for
traces of large random symmetric matrices with independent matrix elements",
\ Bol. Soc. Brasil. Mat. (N.S.) \ \textbf{29}, 1-24 (1998).

\bibitem{So:98} A. Soshnikov, \emph{Level spacings distribution for large
random matrices: Gaussian fluctuations.} Ann. of Math. (2) \textbf{148}
(1998) 573-617

\bibitem{So:00} Soshnikov, A., "The central limit theorem for local linear
statistics in classical compact groups and related combinatorial
identities", \ Ann. Probab. \ \textbf{28}, 1353--1370 (2000).

\bibitem{So:02} A.~Soshnikov, \emph{Gaussian limit for determinantal random
point fields}. Ann. Probab. 30 (2002) 171--187

\bibitem{Sp:87} Spohn, H., "Interacting Brownian particles: a study of
Dyson's model", \ in \emph{Hydrodynamic Behavior and Interacting Particle
Systems} \ edited by G. Papanicolaou (New\ York, Springer, 1987).

\bibitem{Sz:75} Szeg\"{o}, G., \ \emph{Orthogonal Polynomials} \
(Providence, AMS, 1975).

\bibitem{Te:00} Teschl, G., \ \emph{Jacobi Operators and Completely
Integrable Systems} \ (Providence, AMS, 2000).

\bibitem{Wi:02} Wieand, K., "Eigenvalue distribution of random unitary
matrices", \ Prob. Theory Relat. Fields \ \textbf{123}, 202-224 (2002).
\end{thebibliography}
\end{document}